\documentclass[11pt]{article}
\usepackage{geometry}                % See geometry.pdf to learn the layout options. There are lots.
\usepackage{array,multirow}
\usepackage{graphicx}
\usepackage{apalike}
\usepackage{amssymb}
\usepackage{amsmath}
\usepackage{diagbox}
\usepackage{bigints}
\usepackage{epstopdf}
\usepackage{color}
\usepackage{xcolor}
\usepackage{cancel}    % \cancel{ } will strike out it's argument
\usepackage{verbatim}
\usepackage{sidecap}
\usepackage{adjustbox}
\usepackage{mathrsfs}
\usepackage{stmaryrd}
\usepackage{mathtools}
\usepackage{nicefrac}
\usepackage{enumitem}
\usepackage{array}  % for custom column types
\usepackage{ragged2e}  % improves justification tools
\usepackage{multirow}
\usepackage[normalem]{ulem}
\usepackage{tikz}

\newcolumntype{C}[1]{>{\centering\arraybackslash}m{#1}}
\newcolumntype{R}[1]{>{\raggedleft\arraybackslash}m{#1}}

\graphicspath{{figures/}}
\def \figi#1
{\adjustbox{valign=m,vspace=0.7pt}{\includegraphics[scale=0.15]{#1}}}
\usepackage{float}
\usepackage{subcaption}
\usepackage[english]{babel}
\usepackage{blindtext}
\usepackage{authblk}

\title{\textbf{Momentum Equation-Based Regularization and Image Registration for Two-Dimensional Ultrasound Elasticity Imaging}}
\author[1]{Olalekan A. Babaniyi}
\author[1]{Rebecca Rodrigues}
\author[2]{Michael S. Richards}
\affil[1]{School of Mathematics and Statistics, Rochester Institute of Technology, Rochester, NY, USA}
\affil[2]{Department of Biomedical Engineering, Rochester Institute of Technology, Rochester, NY, USA.}

\date{\today}                                           % Activate to display a given date or no date

\begin{document}
\maketitle
%\newpage

%\tableofcontents

%\newpage
%\section{}
%\subsection{}

%---------------------------------------------------------
%  Toggle to label equations 
%
% \def \llabel#1 {\label{#1} }
%\def \llabel#1 {\label{#1} \qquad{#1} }
\def\mean#1{\left< #1 \right>}
\newcommand\myeq{\stackrel{\mathclap{\normalfont(\ref{eq:bvp16})}}{=}}
\newcommand{\ds}{\displaystyle}

%
%---------------------------------------------------------
%---------------------------------------------------------
%  Mark Paul's corrections
%
% \def \llabel#1 {\label{#1} }
%\def \padd#1 {\textcolor{red}{#1}}
\def \pdel#1 {\textcolor{brown}{\sout{#1}}}
\newcommand{\mik}[1]{\textcolor{red}{#1}}
\newcommand{\lek}[1]{\textcolor{blue}{#1}}
\newcommand{\Plots}{1}      %%%%%  PLOTS: 1=ON  0=OFF  %%%%%
%
%---------------------------------------------------------
\newcommand{\bs}[1]{\ensuremath{\boldsymbol{#1}}}

\newcommand{\absdiv}[1]{\par\addvspace{.25\baselineskip}\noindent\textbf{#1}\enspace\ignorespaces}

\def\beq{\begin{equation}}
\def\eeq{\end{equation}}
\def\bea{\begin{eqnarray}}
\def\eea{\end{eqnarray}}
\def\bsig{{\mbox{\boldmath$\sigma$}}}
\def\lam{{\mbox{$\lambda$}}}
\def\tlam{{\mbox{$\widetilde{\lambda}$}}}
\def\tu{{\mbox{$\widetilde{u}$}}}
\def\gam{{\mbox{$\gamma$}}}
\def\bep{{\mbox{\boldmath$\epsilon$}}}
\def\ep{{\mbox{$\epsilon$}}}
\def\kap{{\mbox{$\kappa$}}}
\def\q{{\mbox{\boldmath$q$}}}
\def\r{{\mathbf{r}}}
\def\I{{\mathbf{I}}}
\def\bu{\bs{u}}
\def\buh{\bar{\bs{u}}}
\def\bg{\bs{g}}
\def\bw{\bs{w}}
\def\bx{\bs{x}}
\def\e{{\mbox{\boldmath$e$}}}
\def\K{{\mbox{\boldmath$K$}}}
\def\f{{\mbox{\boldmath$f$}}}
\def\matlab{{\sc{matlab}}}
\def\spreme{{\sc{spreme}}}
\def\tr{\mbox{tr}}
\def\bt{{\mbox{\boldmath$t$}}}
\def\A{\bs{A}}
\def\G{\Gamma}
\def\R{\mathcal{R}}
\def\oh{\omega_{_h}}
\def\uh{u_{_h}}
\def\tk{\tau_{_K}}
\def\ao{\alpha_{\Omega}}
\def\no{\psi_{\Omega}}
\def\spreme{{\sc{spreme}}}
%\def\cC{\mathcal{C}}

%
%---------------------------------------------------------
\begin{abstract}
\absdiv{Objective:} Evaluate and compare multiple mechanics-based and traditional regularization strategies within a variational image registration framework for quasi-static ultrasound elastography.

\absdiv{Methods:}We reformulate a previously proposed momentum-equation-based post-processing method (SPREME) as a regularization term directly integrated into an image registration energy functional. Four regularization types are implemented and compared: a strain magnitude ($\R_\epsilon$), a strain magnitude with incompressibility constraint ($R_{\epsilon i}$), and a momentum-based regularization under plane strain ($R_{P\epsilon}$) and plane stress ($R_{P\sigma}$) assumptions. Each is evaluated in a variational framework solved via Gauss-Newton optimization.

\absdiv{Data:}Registration performance is assessed using synthetic ultrasound image sequences generated from 2D and 3D finite element simulations, as well as experimental phantom data. Comparisons are based on displacement and strain field errors, strain contrast, and contrast-to-noise ratio (CNR).

\absdiv{Results:}Momentum-based regularization, particularly under plane stress assumptions ($R_{P\sigma}$), achieved the lowest strain errors and highest strain contrast across both single-frame and accumulated measurements, even when the underlying tissue deformation violated 2D assumptions. In contrast, strain magnitude regularization with an incompressibility constraint ($R_{\epsilon i}$) produced unstable results in 3D and accumulated displacement scenarios.

\absdiv{Conclusions:}Mechanics-based regularization that incorporates momentum conservation outperforms strain-based techniques in elastographic image registration, particularly when applied directly in the optimization framework. This approach improves robustness to noise and model mismatch, offering a promising direction for future displacement-based inverse imaging methods.
\end{abstract}

\section{Introduction}
Ultrasound elastography (USE), or elasticity imaging, is a technique used to query mechanical states or mechanical properties of soft tissues by extracting information from tissue motion measured from image sequences of deforming tissues \cite{ophir1991elastography,garra2015elastography,shiina2015wfumb,sigrist2017ultrasound}. Quasi-static USE measures tissue motion induced by slowly applying forces externally (e.g., sonographer applied transducer compression) or internally induced forces (e.g., cardiac pressure). Accurate assessment of the deformation field is a critical step in all USE techniques and can have a dramatic effect on subsequent discrimination of the imaged tissue. Traditionally, computational algorithms used to quantify deformations in quasi-static elastography can be loosely classified as block matching methods, energy function optimization-based methods, and deep learning approaches \cite{10551279}. Thorough reviews of many of these approaches can be found in \cite{10551279}, \cite{jiang2018ultrasonic}, and \cite{li2022deep}. 

Block matching algorithms are computationally efficient but suffer from high noise, particularly in estimating lateral displacements or in regions of high strain \cite{10551279}. They can also experience low resolution due to large correlation windows or when using smoothing algorithms designed to compensate for noisy measurements \cite{10551279}. Optimization-based methods, on the other hand, tend to be slower due to the iterative search for displacement measurements \cite{10551279}. To improve measurement accuracy, these techniques will typically incorporate regularization through additional terms in energy functionals that penalize large spatial variations in displacement fields \cite{rivaz2010real,richards2013non}. Many of the regularization strategies implemented for USE penalize first- or second-order spatial derivatives of the displacement field \cite{sumi2008regularization,ashikuzzaman2022second,pellot2004ultrasound}. However, regularization types vary and can significantly bias displacement measurements depending on their form. Physics-inspired regularization terms tend to bias the algorithm output toward a displacement estimate that conforms to the principles or equations that govern deformable bodies and can be incorporated into both deep learning-based and optimization-based approaches \cite{10551279}. Deep learning-based techniques are relatively new to the field of elastography but can potentially merge the accuracy of optimization-based techniques with the computational efficiency of block matching algorithms \cite{li2022deep}. One of the major challenges with deep learning approaches is the need for large training data sets, which, for supervised training, require known solutions obtained from other approaches. As in optimization-based methods, the training of deep learning models involves the minimization of an energy function over the space of a large data set, which often includes a data mismatch term and a regularization term.

Several groups have used physics-inspired regularization, in an optimization framework, to improve displacement estimates in quasi-static USE \cite{poree2015noninvasive,ashikuzzaman2023exploiting,pan2014regularization,mix2017detecting}. In all of these techniques, an assumption of 2D tissue incompressibility was used, via a penalty or a constraint, to ensure a specific mathematical relationship is enforced relating the axial strains to the lateral strains. Gokhale {\it et al.} combined registration and reconstruction, simultaneously searching for a displacement field and a shear modulus distribution for 2D USE \cite{gokhale2004simultaneous}. In that work, the displacement field was constrained by the equations of motion for a linear elastic, incompressible material deforming in 2D plane strain with a given shear modulus distribution. The shear modulus distribution was then optimized to best register sequential image frames. However, this technique requires {\it a priori} the boundary conditions of the problem, which are unknown. Other groups have used a regularized approach to calculate a physically consistent displacement estimate from initial displacement estimates, measured with an unregularized method, adding an extra step to the measurement work flow \cite{kheirkhah2023novel,duroy2023regularization,guo2015pde,babaniyi_spreme}. Babaniyi {\it et al.} developed a SParse RElaxation of the Momentum Equation (SPREME) method to reconstruct lateral displacements using previously measured axial displacement fields only. This method uses a regularization term to globally enforce, but locally relax, the equations of motion for a linear elastic, incompressible material deforming in plane stress. This technique does not require knowing or searching for the shear modulus distribution, it only assumes that the distribution is well represented
as varying in a piecewise constant manner \cite{babaniyi_spreme}.

Whenever regularization is introduced to an estimation, it imposes prior assumptions about the data (e.g. smoothness or physical consistency) and the solution is biased toward that prior. For elastography using 2D imaging techniques, mechanics-based regularization necessarily biases tissue measurements toward 2D assumptions of tissue deformation, where true tissues are in fact deforming in three dimensions. The goal of this work is to develop a momentum-based regularization that penalizes deformations that are inconsistent with 2D plane strain and plane stress using a novel SPREME method. We also implemented two prior approaches to regularizing displacements in 2D elastography. We then compare and contrast these model-based SPREME regularization schemes with a strain-based regularization with and without a penalty imposing 2D tissue incompressibility. We demonstrate the accuracy of the technique using simulated US images pairs and finite element (FE) simulated deformations of linear elastic materials deforming 2D plane strain, 2D plane stress and full 3D deformations. We also show that the SPREME regularizations can be used to improve accumulated estimates of displacement fields over large strains, measured from many sequential frames. Lastly, we present an example measurement of displacement and strain fields from a 3D tissue mimicking phantom designed for USE validation. 

\section{Image Registration Formulation}
The measurement algorithm developed here is an optimization problem, where the goal is to find a displacement field $\bu(\bx)$, within the measurement domain $\Omega \subset \mathbb{R}^2$, that minimizes the following objective function:
\beq
\pi[\bu_k(\bx)] = \psi[\bu_k(\bx)] + \mathcal{R}[\bu_k(\bx)] \qquad k = 1, 2, \ldots N-1, \label{eq:objf}
\eeq
where $N \geq 2$ is the number of frames, $\bu_k(\bx)$ is the incremental frame-to-frame displacement, and $\psi[\bu_k(\bx)]$ is the image matching term defined as
\beq
\psi[\bu_k(\bx)] = \frac{1}{2}\int_{\Omega} (I_k(\bx + \bs{S}_k(\bx)) - I_{k+1}(\bx + \bs{S}_k(\bx) + \bu_k(\bx)))^2 \ d\Omega, \label{eq:IR}
\eeq
where $\bs{S}_k(\bx) = \sum_{i=0}^{k-1} \bu_i(\bx)$ is the accumulated displacement up to frame $k-1$ and $\bu_0(\bx)=\bs{0}$. This term quantifies the difference (or mismatch), in the measurement domain $\Omega$, between an undeformed reference image, $I_1(\bx + \bs{s}_k(\bx))$, and a target image $I_2(\bx + \bs{s}_k(\bx) + \bu_k(\bx))$ that has been non-linearly deformed by the displacement field measurement $\bu_k(\bx)$. For the remainder of the derivation, we focus on the special case where $N=2, \bs{S}_1(\bx) = \bs{0}$, and we set $\bu_k = \bu$ to simplify the notation. It should be noted that for the general case where $N>2$, the optimization problem described in the special case is solved repeatedly and sequentially. Furthermore, all the $N-1$ accumulated displacement fields will share a common reference coordinate frame. The regularization term $\R[\bu(\bx)]$ is used to deal with the ill-posedness of the registration problem. Typically, there are multiple $\bu(\bx)$ that can be used to match images $I_1$ and $I_2$. Furthermore, small changes in the images can result in very different displacement fields.

The regularization term can be used to obtain unique and stable solutions by specifying desired features of the solution (e.g., smoothness, small magnitude, etc.). It can also be used to provide additional information about the optimization problem, such as physical constraints that $\bu(\bx)$ need to satisfy. In this work, we consider different types of regularization terms that achieve one or both of these goals, as described in the next section.

\subsection{Regularization Types}
We consider four different types of regularization terms to help reconstruct unique displacement fields with desired properties. The first one is a strain magnitude regularization defined as
\beq
\R_{\epsilon}[\bu(\bx)] = \frac{\alpha_{\epsilon}}{2} \int_{\Omega} \nabla^s \bu(\bx): \nabla^s \bu(\bx) \ d\Omega,
\eeq 
where
\beq
\nabla^s \bu = \frac{1}{2}\left(\nabla \bu + (\nabla \bu)^T\right), \label{eq:strain}
\eeq
and $\bs{B}:\bs{C}$ is the inner product of two tensor fields. This regularization term has previously been used in \cite{richards2013non,richards2009quantitative}, and penalizes large strains, with the constant $\alpha_{\epsilon}$ regulating how strongly it is weighted relative to the image matching term.

The second regularization type we implement is one that penalizes large strains and enforces an incompressibility constraint via a penalty term \cite{mix2017detecting}. We refer to this as the strain regularization with incompressibility constraint ($\R_{\epsilon i}$) and define it as
\beq
\R_{\epsilon  i}[\bu(\bx)] = \frac{\alpha_{\epsilon i}}{2} \int_{\Omega} \nabla^s \bu(\bx): \nabla^s \bu(\bx) \ d\Omega + \frac{\alpha_i}{2} \int_{\Omega} (\nabla \cdot \bu(\bx))^2 \ d\Omega, \label{eq:saicr}
\eeq
where $\alpha_{\epsilon i}, \alpha_i$ are constants penalizing strains and enforcing the incompressibility constraint, respectively. This regularization term has previously been applied to the image registration problem in \cite{richards2013non,richards2009quantitative}, and used to reconstruct precise estimates of the full displacement field in some situations. For 2D USE, it is an appropriate constraint when the material being imaged is deforming in plane strain and incompressible. 

The final two regularization terms we consider are based on enforcing the conservation of momentum equations on the displacement field. The derivation of the momentum-based regularization terms developed here is similar to the post-processing algorithm developed in \cite{babaniyi_spreme} with several exceptions, as will be noted in this section. It is assumed that the 2D vector field within the region $\Omega$ is well approximated by a plane stress or a plane strain approximation of a linear elastic, isotropic, incompressible or nearly incompressible material. Furthermore, it is assumed that body forces are negligible, and that the deformation is small and quasi-static, thus giving the following equilibrium equation:
\beq
\nabla \cdot \bsig =\bs{0} \qquad \mbox{in } \Omega,
\eeq
where the stress tensor, $\bsig(\bx)$, is defined as
\beq
\bsig(\bx) = \mu(\bx) \bar{\lambda}\nabla \cdot \bu(\bx) \I + 2 \mu(\bx) \nabla^s \bu(\bx), \label{eq:cEq}
\eeq  
and $\mu(\bx)$ is the shear modulus distribution. For an incompressible material (i.e., Poisson's ratio $\nu = 0.5$) that is deforming in 
plane stress, the parameter $\bar{\lambda} = 2$. For a nearly incompressible material that is deforming in plane strain, the parameter $\bar{\lambda} = \frac{2\nu}{1-2\nu}$,
noting that as $\nu \rightarrow 1/2$, $\bar{\lambda} \rightarrow \infty$. Thus, Eq.\@ (\ref{eq:cEq}) represents a generalized model for both incompressible plane stress and nearly incompressible plane strain depending on the choice of $\bar{\lambda}$, either $\bar{\lambda} = 2$ or $\bar{\lambda} \gg 2$, respectively. The equilibrium equation can then be written as
\beq
\nabla \cdot (\mu \bs{A}) = \bs{0} \qquad \mbox{in } \Omega, \label{eq:meq}
\eeq
where,
\beq
\bs{A}[\bu] = \bar{\lambda} \nabla \cdot \bu(\bx) \I+ 2 \nabla^s\bu(\bx).
\eeq
It was shown in \cite{babaniyi_spreme} that if $\mu(\bx)$ is a piecewise constant function, then Eq.\@ (\ref{eq:meq}) reduces to
\beq
\nabla \cdot \bs{A}[\bu] = \bs{0} \qquad a.e. \qquad \mbox{in } \Omega, \label{eq:pwcmeq}
\eeq
where a.e. means \textit{almost everywhere} in the domain or everywhere where there is a positive measure (i.e., excluding curves and points). Since $\nabla \cdot \bs{A}$ is sparse, it can be enforced using a Total Variation (TV) regularization as shown below:
\beq
\R_{Pm}[\bu(\bx)] = \alpha_{Pm} \int_{\Omega} ||\nabla \cdot \bs{A}(\bx)||_2 \ d\Omega, \label{eq:mer}
\eeq
where $||\cdot||_2$ is the Euclidean norm of a vector field in $\mathbb{R}^2$, and $\alpha_{Pm}$ is the regularization parameter used to regulate how strongly Eq.\@ (\ref{eq:mer}) is enforced relative to the image matching term. The derivative of the absolute value function is undefined at the origin; hence an approximation is needed if one wants to use an optimization algorithm that requires derivatives to minimize Eq.\@ (\ref{eq:objf}). In \cite{babaniyi_spreme}, the authors use a quadratic approximation of the absolute value. Here, we use the following approximation:
\beq
\R_{Pm}^{\delta}[\bu(\bx)] \approx \alpha_{Pm} \int_{\Omega} \sqrt{(\nabla \cdot \A[\bu]) \cdot (\nabla \cdot \A[\bu]) + \delta} \ d\Omega, \label{eq:mreg}
\eeq 
where $\delta$ is a small positive constant used to make $\R_{Pm}^{\delta}[\bu(\bx)]$ differentiable when $\nabla \cdot \A[\bu] = \bs{0}$.

\subsection{Optimization}
We solve the optimization problem Eq.\@ (\ref{eq:objf}) using Newton's method.  Specifically, we use a Gauss-Newton approximation for the image matching term, and an approximation for the second derivatives of the momentum regularization term Eq.\@ (\ref{eq:mreg}) as will be explained later on in this section. Given an initial guess $\bu_1$,  $\bu_j, j = 1,2,\ldots,$ is updated in each iteration $j$ by:
\beq
\bu_{j+1} = \bu_j + \delta\bu_j,
\eeq
where $\delta\bu_j$ is the search direction. This search direction is determined by solving the linear system:
\beq
(D^2_{\bu \bu} \pi[\bu_j]\delta\bu_j, \bw) = -(D_{\bu}\pi[\bu_j], \bw) \qquad \text{for all}\ \bw,
\eeq
where $D_{\bu}\pi$ and $D^2_{\bu \bu} \pi$ are the first and second derivatives of the functional Eq.\@ (\ref{eq:objf}), respectively, and $(a, b) = \int_{\Omega} ab \ d\Omega$ is the $L^2(\Omega)$ inner product. The derivatives can be expressed in terms of the data match and regularization terms as
\begin{subequations}
\begin{alignat}{2}
(D_{\bu}\pi[\bu], \bw) &= (D_{\bu}\psi[\bu], \bw) + (D_{\bu}R[\bu], \bw) \qquad &&\text{for all}\ \bw, \\
(D^2_{\bu\bu}\pi[\bu]\delta\bu, \bw) &= (D^2_{\bu\bu}\psi[\bu]\delta\bu, \bw) + (D^2_{\bu\bu}R[\bu]\delta\bu, \bw) \qquad &&\text{for all}\ \bw, \delta\bu,
\end{alignat}
\end{subequations}
where $\bs{w}$ and $\delta \bs{u}$ are arbitrary admissible variations of $\bu$ \cite{hughes2012finite}, and $(\bu, \delta \bu) = (\bu_j, \delta \bu_j)$. We drop $j's$ in most of the rest of this section to simplify the notation. The first and second derivatives of $\psi[\bu]$ are:
\begin{subequations}
\label{eq:imds}
\begin{alignat}{2}
(D_{\bu}\psi[\bu], \bw) &= \int_{\Omega_I} \bw \cdot \nabla I_2(\bx + \bu(\bx))(I_2(\bx + \bu(\bx)) - I_1(\bx)) \ d\Omega\qquad &&\forall \bw,  \\
(D^2_{\bu\bu}\psi[\bu]\delta\bu, \bw) &= \int_{\Omega_I} \bw \cdot (\nabla I_2(\bx + \bu(\bx))\otimes \nabla I_2(\bx + \bu(\bx))) \delta \bu \ d\Omega \qquad &&\forall \bw, \delta \bu. \label{eq:imh}
\end{alignat}
\end{subequations}
We used the Gauss-Newton approximation and dropped the second derivatives of the images in Eq.\@ (\ref{eq:imh}). To improve computational efficiency, it is further assumed that when $\bu_j (\bx)$ is sufficiently close to the optimal displacement field (that is, at the minimum of eqn (\ref{eq:objf})) then $\nabla I_2 (\bx +\bu_j(\bx) ) \approx \nabla I_1(\bx)$. This assumption eliminates the need to interpolate $\nabla I_2(\bx +\bu_j(\bx))$ at every iteration. However, it is still necessary to interpolate $I_2(\bx +\bu_j(\bx))$ \cite{rivaz2010real}. Using this approximation in Eq.\@ (\ref{eq:imds}) gives:
\begin{subequations}
\label{eq:wfnew}
\begin{alignat}{2}
(D_{\bu}\psi[\bu], \bw) &= \int_{\Omega_I} \bw \cdot \nabla I_1(\bx)(I_2(\bx + \bu(\bx)) - I_1(\bx)) \ d\Omega\qquad &&\forall  \bw,  \\
(D^2_{\bu\bu}\psi[\bu]\delta\bu, \bw) &= \int_{\Omega_I} \bw \cdot (\nabla I_1(\bx)\otimes \nabla I_1(\bx)) \delta \bu \ d\Omega \qquad &&\forall \bw, \delta \bu.
\end{alignat}
\end{subequations}

Next, we provide the derivatives of the different regularization types. We note that when $\alpha_i = 0$ in (\ref{eq:saicr}), $\R_{\epsilon i}[\bu]$ reduces to $\R_{\epsilon}[\bu]$. Hence, we only provide the first and second derivatives for $\R_{\epsilon i}[\bu]$ below:
\begin{subequations}
\begin{alignat}{2}
(D_{\bu}\R_{\epsilon i}[\bu], \bw) &= \alpha_{\epsilon i}\int_{\Omega} \nabla^s \bs{w} : \nabla^s \bu \ d\Omega + \alpha_i\int_{\Omega} (\nabla \cdot \bw) (\nabla \cdot \bu) \ d\Omega &&\hspace{2mm} \forall \bw, \\
(D^2_{\bu\bu}\R_{\epsilon i}[\bu]\delta\bu, \bw) &= \alpha_{\epsilon i}\int_{\Omega} \nabla^s \bs{w} : \nabla^s \delta\bu \ d\Omega + \alpha_i\int_{\Omega} (\nabla \cdot \bw) (\nabla \cdot \delta\bu) \ d\Omega &&\hspace{2mm} \forall \bw, \delta \bu.
\end{alignat}
\end{subequations}

The first and second derivatives of the momentum-based regularization (\ref{eq:mreg}) are:
\begin{subequations}
\begin{alignat}{2}
(D_{\bu}\R^{\delta}_{Pm}[\bu], \bw) &= \alpha_{Pm} \int_{\Omega} \frac{(\nabla \cdot \A[\bw])\cdot (\nabla \cdot \A[\bu])}{\sqrt{(\nabla \cdot \A[\bu]) \cdot (\nabla \cdot \A[\bu]) + \delta}} \ d\Omega \qquad \forall \bw, \label{eq:gmr}\\
(D^2_{\bu\bu}\R^{\delta}_{Pm}[\bu]\delta\bu, \bw) &= \alpha_{Pm} \int_{\Omega} \frac{1}{\sqrt{\|\bg(\bu)\|^2_2 + \delta}}\left[\left(\I - \frac{\bg(\bu)\otimes \bg(\bu)}{\|\bg(\bu)\|^2_2 + \delta} \right) \nabla \cdot \bg(\bw) \right]\cdot \bg(\delta \bu) \ d\Omega \notag,
\end{alignat}
\end{subequations}
where $\bg(\bu) = \nabla \cdot \A[\bu]$ is used to express $D^2_{\bu\bu}\R^{\delta}_{Pm}$ in compact form. The nonlinear term $ \ds \bg(\bu)\otimes \bg(\bu) (\|\bg(\bu)\|^2_2 + \delta)^{-1}$ in the second derivative has been shown to cause convergence problems for Newton's method \cite{chan1995continuation,vogel1996iterative}. We therefore remove this term, and use the following approximation for the second derivative \cite{chan1995continuation}:
\beq
(D^2_{\bu\bu}\R^{\delta}_{Pm}[\bu]\delta\bu, \bw) \approx \alpha_{Pm} \int_{\Omega} \frac{(\nabla \cdot \A[\bw])\cdot (\nabla \cdot \A[\delta\bu])}{\sqrt{(\nabla \cdot \A[\bu]) \cdot (\nabla \cdot \A[\bu]) + \delta}} \ d\Omega \qquad \forall \bw, \delta \bu. \label{eq:hmr}
\eeq

For all regularization types used here, the displacement fields are discretized using bilinear quadrilateral FE shape functions. In Eq.\@ (\ref{eq:gmr}) and Eq.\@ (\ref{eq:hmr}), $\A[\bu]$ has an explicit second spatial derivative of $\bu(\bx)$, which is incompatible with the FE discretizations. In \cite{babaniyi_spreme}, they mitigated this problem by defining their optimization problem in terms of $\bu(\bx)$ and strains $ \bep(\bx)$, separately, and then expressing $\A$ in terms of $\bep$ (i.e. $\A[\bep]= \bar{\lambda} tr(\epsilon) \I+ 2\epsilon$). In that work, $\bep(\bx)$ and $\bu(\bx)$ were coupled by means of an additional compatibility term, adding to their computational expense. In this work, we chose to mitigate this issue by utilizing linear interpolating functions for our searched displacement fields and assuming $\A[\bu]$ is approximately piecewise constant\footnote{The displacement fields in this work are discretized with bilinear, quadrilateral shape functions and thus $\A$ within elements is only approximately piecewise constant (see Appendix Section \ref{apss:mbr}).}. We follow the approach presented in \cite{richards:thesis,richards2009quantitative} to derive expressions for Eq.\@ (\ref{eq:gmr}) and Eq.\@ (\ref{eq:hmr}) for the discretized problem. It should be noted here that the method of solving and discretizing the strain-based regularizations are not quite the same as the SPREME regularizations. Due to the unique way in which we chose to reduce Eq.\@ (\ref{eq:mer}), the optimization performed here is specifically dependent on the discretization of $\bu(\bx)$ used here. A detailed description of the derivation and assumptions is given in the Appendix Section \ref{apss:mbr}.

For optimization algorithms, it is necessary to use a global search algorithm, such as a block matching method, to establish an initial guess of the displacement field $\bu_1$ \cite{kheirkhah2023novel,poree2015noninvasive,richards2013non}. Here, a simple 2D normalized cross-correlation block matching algorithm was used to find an estimate of $\bu_1$ \cite{richards2013non}. A measurement window size of 5 mm axially and 9 mm laterally was used, with 25\% and 40\% overlap, respectively. Initial displacements were only estimated to the nearest integer value of pixels. A 5$\times$5 pixel median filter was applied to each image of the displacement component to reduce outliers in the estimate. The measured displacement vector field was then linearly interpolated onto the FE mesh and used as the initial guess of the displacement estimate in the optimization.

The regularization parameters ($\alpha_{\epsilon}$, $\alpha_{\epsilon i}$,  $\alpha_{P\epsilon}$, and $\alpha_{P\sigma}$) were chosen to be the values that minimized the error in the total strain measurement, for a single-frame registration, when the frame-to-frame strain was approximately $0.8\%$, which approximately corresponds to the maximum frame-to-frame strain in the simulated data generated in this work. The values were $\alpha_{\epsilon}=24.8$, $\alpha_{\epsilon i}=1.89e4$, $\alpha_{P\epsilon}=1.14e\text{-}4$, and $\alpha_{P\sigma}=2.33e\text{-}4$. See Appendix Section \ref{apss:par} for details of methodology for choosing regularization parameters. Furthermore, for the strain regularization with incompressibility constraint ($\R_{\epsilon i}$) $\alpha_{i}=100 \times \alpha_{\epsilon i}$. For the approximation of momentum-based regularization, we used $\delta = 1e\text{-}8$. For the plane strain regularization ($\R_{P\epsilon}$) $\bar{\lambda}=9$, which corresponds to $\nu=0.45$, and for the plane stress regularization ($\R_{P\sigma}$) $\bar{\lambda}=2$.

The frame-to-frame optimizations were stopped at $j = 20$ iterations or when $\|\delta \bu_j\|/\|\bu_j\| < 0.001$, whichever occurred first. The registration algorithm was written and implemented in Matlab (Mathworks, Natick, MA). The 2D finite element meshes were created using Matlab generated input files to an open source FE meshing software (Gmsh) \cite{geuzaine2009gmsh}.

\section{Data Generation}
\subsection{Simulated Data}
\label{ssec:sd}
Computer-simulated radio-frequency (RF) US images were created using a point spread function (PSF) approximation to model US image formation and a distribution of point scatterers within a field. Scatterers movement, to create sequences US frames, was determined by generating known displacement fields using commercial FE software (Comsol Inc., Stockholm, SE). FE simulations were designed to mimic the geometry and position of the image that could be achieved in free-hand US breast imaging, where a localized inclusion lies within a homogeneous background material, imaged with a linear array US transducer (Figure \ref{fig0}(a)). The tissue was modeled as a linear elastic, nearly incompressible material undergoing quasi-static deformation. The geometry of the simulated tissue consisted of a block of material with dimensions 16$\times$16$\times$16 cm and a Young’s modulus of 10 MPa. Within the material block, a spherical inclusion with a diameter of 1.5 cm was centrally located 3 cm from the top of the model, as shown in Figure \ref{fig0}(a). Multiple models with inclusion Young’s moduli of 40 MPa (hard inclusion, HI), 20 MPa (medium hard inclusion, mHI), 10 MPa (homogeneous model, HOM) and 5 MPa (soft inclusion, SI) were generated separately to compare models with varying inclusion magnitudes. The Poisson's ratio of the material was set to $\nu$=0.495. The deformations were induced by fixing the bottom displacement of the model in all directions and compressing the block axially (y direction) within a square region 5$\times$5 cm in the center of the top of the block. A slip boundary condition was applied on this top compressed region. This was intended to mimic transducer compression from the skin’s surface with an attached compression plate using coupling gel between the transducer and the skin surface \cite{varghese2002elastographic}. The remaining surfaces of the block were traction-free. The compression resulted in the total overall strain within the simulated image plane shown in  (Figure \ref{fig0}(b)). Two-dimensional forward simulations were also performed using geometries and boundaries equivalent to the imaged plane of the 3D geometry (HI only) described above. In these cases, the central slices were extracted and modeled as deforming in either plane strain or plane stress. These 2D forward models were created to quantify and compare measurement accuracy in situations where our model assumptions are true rather than an approximation.
\ifnum \Plots>0
\begin{figure}[h]
    %\centering
    \includegraphics[trim={0cm 1.5cm 0cm 4cm},clip,width = \textwidth]{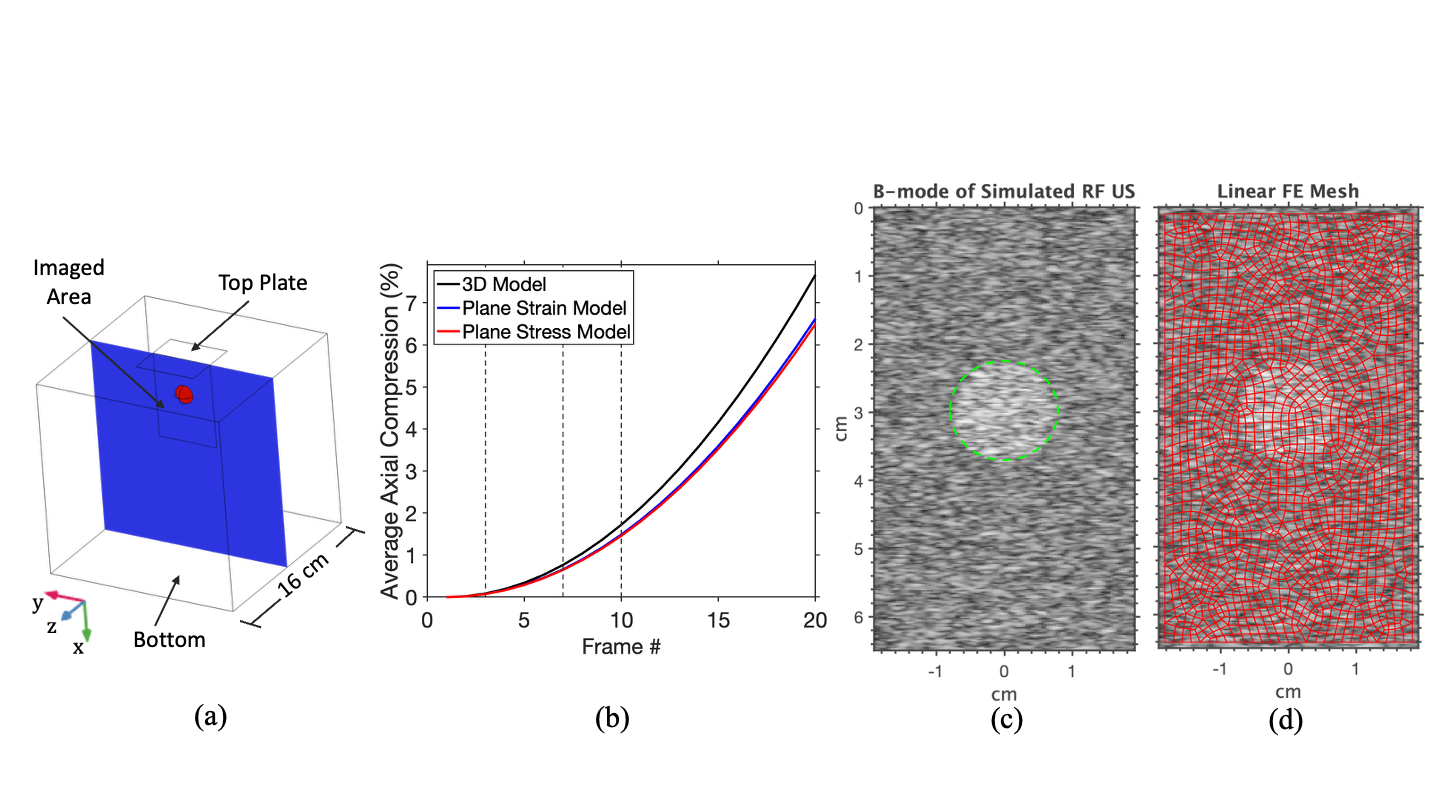}
    \caption{(a) Finite element model geometry showing top plate, used for applying deformation to cube, the central slice used for 2D simulations of plane strain and plane stress deformations (blue), the relative location of the spherical (or circular in 2D) inclusion (red) and the boundaries of the imaged area use to simulate US images. (b) Average axial strain ($\epsilon_{xx}$) within the imaged plane in each frame for the 3D and 2D models. The dotted lines indicate frames 3, 7 and 10 used for dual frame displacement measurements. (c) Example B-mode image create from the envelop of the simulated RF images used in the study. (d) Finite element mesh (red) used for discretizing the measured displacements in the image registration.}
    \label{fig0}
\end{figure}
\fi
Simulated RF US image pairs were created for all simulations using a 2D, spatially invariant PSF. The PSF approximated an US imaging pulse for a linear array transducer with a center frequency of 5.5 MHz, a Gaussian modulated pulse length of 0.43 $\mu$s and a lateral full width at half maximum of 1.4 mm. This PSF was convolved with evenly distributed point scatterers (density 30 mm$^{\text{-}2}$) with an arbitrarily assigned amplitude within the imaged area (4.4 cm laterally by 7.2 cm axially), as shown in Figure \ref{fig0}(a), to create a reference image, $I_1(\bx)$. The scatterers were then displaced according to the simulated displacement fields and convolved with the same PSF to create a target image, $I_2(\bx)$. The scattering amplitudes of the points within the inclusion region were increased (2$\times$) to improve visualization. The images were sampled axially at the equivalent of a 40 MHz sampling frequency (at 1540 m/s, spatial sampling was 0.0192 mm) and laterally with a spacing of 0.3 mm. Normally distributed Gaussian noise was added to the final RF US images such that the images had an approximate signal to noise ratio of SNR = 12 dB. Figure \ref{fig0}(c) shows an example B-mode image generated from the simulated RF linear array image. The displacement field was then interpolated from the COMSOL FE mesh to the scatterer positions and used to model the deformed scatterer locations. The target image, $I_2 (x)$, was then created by convolution of the US-simulated PSF, as above, with the same scatterers at the new locations. The target images had the exact same sample spacing and added noise as the reference image. The image sequences of 20 total frames were created for each forward model by scaling the displacement fields so that the average axial compressive strain from frame to frame within the image domain was approximately $0.4\%$, the maximum was $0.8\%$, and the total strain over the 20 frame sequence was $\sim7\%$. Figure \ref{fig0}(b) shows the average axial compressive strain, within the imaged plane, for all frames in the sequence for each of the three forward models. Figure \ref{fig0}(d) shows the mesh used for all the simulated data registrations performed in this work. All displacement field measurements from the simulated data used the same mesh, and simulated image sequences so that a comparison across regularization types could be made.

\subsection{Phantom Experiments}
To test our method on real-world data, we captured images from a tissue-mimicking phantom designed for USE validation (Elasticity QA Phantom, Model 049, CIRS Inc., Norfolk, VA, USA). We imaged a 20mm diameter spherical inclusion with an elastic modulus of $80 \pm 12$ kPa embedded in a background material of $25 \pm 6$ kPa. This inclusion had an expected modulus contrast of approximately 3.2. SonixTouch Research ultrasound system (Ultrasonix Medical Corporation, BC, Canada) was used to acquire data using an L14-5/38 linear array transducer operating at a center frequency of 10 MHz. RF data were collected at a sampling rate of 40 MHz, a lateral spacing of 0.3 mm, and a frame rate of 123 Hz. Each frame consisted of 128 scan lines and 2864 axial samples per line, covering a 55 mm × 38 mm axial by lateral image area \cite{tuladhar2025ultrasound}. The image sequences were decimated in time and accumulated in 13 frames, such that the average axial compressive frame-to-frame strain was approximately $0.35\%$ and the total axial compressive strain was approximately $\sim4\%$. An example B mode and the FE mesh generated for the phantom USE images are shown in Figures \ref{fig:phantom_og}(a) and \ref{fig:phantom_og}(b), respectively. All displacement field measurements from the phantom data used the same mesh and image sequence to compare across regularization types.

\newcommand{\cw}{1cm}
\newcommand{\ctw}{0.95cm}

\section{Results}
To evaluate the impact of regularization strategies on image registration performance, we performed a comparative analysis across the four types of regularization (i.e. $\R_\epsilon$, $\R_{\epsilon i}$, $\R_{P\epsilon}$, and $\R_{P\sigma}$). The performance of each approach was assessed using a percent displacement error estimate as follows \cite{babaniyi_spreme}:
\beq
\text{Disp. Error}_i = \frac{\|u^t_i-u^m_i\|}{\|u^t_i\|}=\frac{\sqrt{\int_\Omega{(u^t_i-u^m_i)^2}\,d\Omega}}{\sqrt{\int_\Omega{(u^t_i)^2}\,d\Omega}}
\eeq
for the $i^{th}$ ($x$ or $y$) component of the displacement vector field, or
\beq
\text{Total Disp. Error} = \frac{\|\bu^t-\bu^m\|}{\|\bu^t\|}=\frac{\sqrt{\int_\Omega{(\bu^t-\bu^m)\cdot(\bu^t-\bu^m)} \,d\Omega}}{\sqrt{\int_\Omega{(\bu^t\cdot\bu^t)}\,d\Omega}},
\eeq
for the total displacement field error. Here, $u^t_i$ and $u^m_i$ are the $i^{th}$ component of the true, forward simulated and measured displacement fields, respectively, and $\bu^t$ and $\bu^m$ are the total true, forward simulated and measured displacement fields. Similarly, the percent strain error was calculated as \cite{babaniyi_spreme}
\beq
\text{Strain Error}_{ij} = \frac{\|\epsilon^t_{ij}-\epsilon^m_{ij}\|}{\|\epsilon^t_{ij}\|}
\eeq
for the ${ij}^{th}$ ($xx$, $yy$, or $xy$) component of the strain tensor field, or
\beq
\text{Total Strain Error} = \frac{\|\bep^t-\bep^m\|}{\|\bep^t\|}=\frac{\sqrt{\int_\Omega{(\bep^t-\bep^m):(\bep^t-\bep^m)} \,d\Omega}}{\sqrt{\int_\Omega{(\bep^t:\bep^t)}\,d\Omega}},
\eeq
where the strain tensors $\bep^t$ and $\bep^m$ are calculated from the true and measured displacement fields, respectively, using eq. \ref{eq:strain}, where $\bep=\nabla^s \bu$. All error calculations were integrated on the domain of the FE mesh.

In addition to error, the contrast in the axial strain was calculated using the strain ratio ($SR$) defined as \cite{shehata2022qualitative}
\beq
SR = \frac{\epsilon_A}{\epsilon_B}=\frac{\int_{\Omega_A}{\epsilon^m_{xx}\,d\Omega_A}}{\int_{\Omega_B}{\epsilon^m_{xx}\,d\Omega_B}},
\eeq
where $\epsilon_A$ and $\epsilon_B$ are the average measured axial strain ($\epsilon^m_{xx}$) outside ($A$) and within ($B$) an identified inclusion region of interest (ROI), respectively. Lastly, the elastographic contrast-to-noise ratio ($CNR_e$) was calculated as \cite{bilgen1997predicting}
\beq
CNR_e = \frac{2(\epsilon_A-\epsilon_B)^2}{(\eta^2_A+\eta^2_B)}.
\eeq
Here, $\eta_A$ and $\eta_B$ are the variances of the axial strain measured outside and within the ROI, respectively. The ROI of the simulated images is identified as the green dotted line in the B mode image shown in Figure \ref{fig0}(c) and the ROI of the phantom images is shown with a similar line in Figure \ref{fig:phantom_og}(a). 

\subsection{Single-frame Displacement Measurements}
The following results correspond to a single-frame registration of 2 US frames from the simulated data sequences. That is, registering frame 7 to frame 1. Table \ref{table:strs} shows the calculated values of the displacement and strain errors as well as the contrast and contrast to noise metrics for the hardest inclusion (HI) and all simulated model sequences (3D, plane strain and plain stress). Table \ref{table:Astrs}, found in the Appendix \ref{apss:suppres}, shows the calculated values of the error and contrast metrics for the soft (SI), homogeneous (HOM) and medium hard inclusions (mHI) for the 3D simulated model sequences. Images of the axial ($u_x$) and lateral ($u_y$) displacement and strain images for all regularization types, along with the true, forward simulated displacement fields, are shown in Figure \ref{fig1}. The three components of the corresponding strain images ($\epsilon_{xx}$, $\epsilon_{yy}$, and $\epsilon_{xy}$) are shown in Figure \ref{fig2}. Axial strain images ($\epsilon_{xx}$) recovered for all regularization types for the plane strain and plane stress forward simulations (HI) are shown in Figure \ref{fig3} with the true strain images. Images of the lateral and shear strains recovered from the plane strain and plane stress forward simulations (HI) can be found in the Appendix \ref{apss:suppres} in Figures \ref{fig:Alat} and \ref{fig:Ashear}, respectively. Axial strain images recovered for all regularization types for the soft and medium hard inclusion forward simulations (3D model only) are shown in Figure \ref{fig4}.
\begin{table}[h]
\begin{center}
\begin{tabular}{|R{0.75cm}||C{\cw}|C{\cw}|C{\cw}|C{\cw}|C{\cw}|C{\cw}|C{\cw}|C{0.7cm}|C{0.95cm}|} 
 \hline
 \multicolumn{8}{|c|}{\bf{Measurement Error (Single-frame Meas. 1 to 7, HI)}}&\multicolumn{2}{|c|}{\bf{Contrast}}\\
 \hline
 \multicolumn{1}{|c||}{Model/Reg.}&\multicolumn{2}{|c|}{Disp. Comp.} &Total Disp.&\multicolumn{3}{|c|}{Strain Component} &Total Strain&\multicolumn{2}{|c|}{$\epsilon_{xx}$ Only}\\
 \hline\hline
 \multicolumn{1}{|l||}{\bf{3D}}  & $u_{x}$ & $u_{y}$  & $\mathbf{u}$ & $\epsilon_{xx}$ & $\epsilon_{yy}$ & $\epsilon_{xy}$  & $\epsilon$&$SR$&$CNR_e$\\
 \hline\hline
$R_\epsilon$ & $0.2\%$ & $20.5\%$ & $3.0\%$ & $14.9\%$ & $45.0\%$ & $60.4\%$ & $26.2\%$ & $1.68$ & $2.72$ \\
\hline
$R_{\epsilon i}$ & $0.3\%$ & $97.0\%$ & $14.0\%$ & $6.0\%$ & $92.1\%$ & $75.6\%$ & $43.8\%$ & $1.54$ & $2.54$ \\
\hline
$R_{P\epsilon}$ & $0.1\%$ & $17.1\%$ & $2.5\%$ & $6.2\%$ & $20.8\%$ & $55.6\%$ & $14.3\%$ & $1.67$ & $3.53$ \\
\hline
$R_{P\sigma}$ & $0.1\%$ & $9.0\%$ & $1.3\%$ & $4.9\%$ & $16.5\%$ & $33.6\%$ & $10.3\%$ & $1.67$ & $3.68$ \\
 \hline\hline
 \multicolumn{1}{|l||}{\bf{P. Strain}} & $u_{x}$ & $u_{y}$  & $\mathbf{u}$ & $\epsilon_{xx}$ & $\epsilon_{yy}$ & $\epsilon_{xy}$  & $\epsilon$&$SR$&$CNR_e$\\
 \hline\hline
 $R_\epsilon$ & $0.3\%$ & $11.8\%$ & $3.3\%$ & $18.7\%$ & $36.7\%$ & $77.1\%$ & $31.1\%$ & $2.04$ & $5.19$ \\
\hline
$R_{\epsilon i}$ & $0.3\%$ & $3.0\%$ & $0.9\%$ & $5.9\%$ & $5.7\%$ & $24.5\%$ & $6.9\%$ & $1.90$ & $7.21$ \\
\hline
$R_{P\epsilon}$ & $0.1\%$ & $4.3\%$ & $1.2\%$ & $3.8\%$ & $5.5\%$ & $19.6\%$ & $5.6\%$ & $1.94$ & $8.93$ \\
\hline
$R_{P\sigma}$ & $0.1\%$ & $6.7\%$ & $1.9\%$ & $5.2\%$ & $11.7\%$ & $37.2\%$ & $10.7\%$ & $2.02$ & $8.93$ \\
 \hline\hline
 \multicolumn{1}{|l||}{\bf{P. Stress}} & $u_{x}$ & $u_{y}$  & $\mathbf{u}$ & $\epsilon_{xx}$ & $\epsilon_{yy}$ & $\epsilon_{xy}$  & $\epsilon$&$SR$&$CNR_e$\\
 \hline\hline
 $R_\epsilon$ & $0.3\%$ & $21.0\%$ & $2.6\%$ & $16.4\%$ & $45.1\%$ & $66.6\%$ & $25.7\%$ & $2.49$ & $11.5$ \\
\hline
$R_{\epsilon i}$ & $0.3\%$ & $136\%$ & $17.1\%$ & $7.5\%$ & $133\%$ & $98.8\%$ & $57.3\%$ & $2.2$ & $14.2$ \\
\hline
$R_{P\epsilon}$ & $0.2\%$ & $16.8\%$ & $2.1\%$ & $6.1\%$ & $23.2\%$ & $60.4\%$ & $14.1\%$ & $2.47$ & $21.9$ \\
\hline
$R_{P\sigma}$ & $0.1\%$ & $9.9\%$ & $1.3\%$ & $3.5\%$ & $14.6\%$ & $30.9\%$ & $8.2\%$ & $2.45$ & $24.8$ \\
 \hline
\end{tabular}
\caption{Percent displacement and strain error calculated for a single-frame measurement (frames 1 to 7) and the hard inclusion (HI) for all forward simulations and shown for all regularization types. Results are given for individual strain tensor components and total strain. Contrast recovery is also reported as $SR$ and $CNR_e$. The strain ratio values calculated from the true strains for the 3D, plane strain and plane stress forward models are $SR = 1.65$, $2.01$, and $2.45$ respectively. The corresponding contrast to noise ratios are $CNR_e = 3.77$, $9.20$, and $24.0$, respectively.}
\label{table:strs}
\end{center}
\end{table}
\ifnum \Plots>0
\begin{figure}[h]
    %\centering
    \begin{tikzpicture}
        \node (image) at (0,0) {\includegraphics[trim={0cm 0 4cm 0},clip,width = \textwidth]{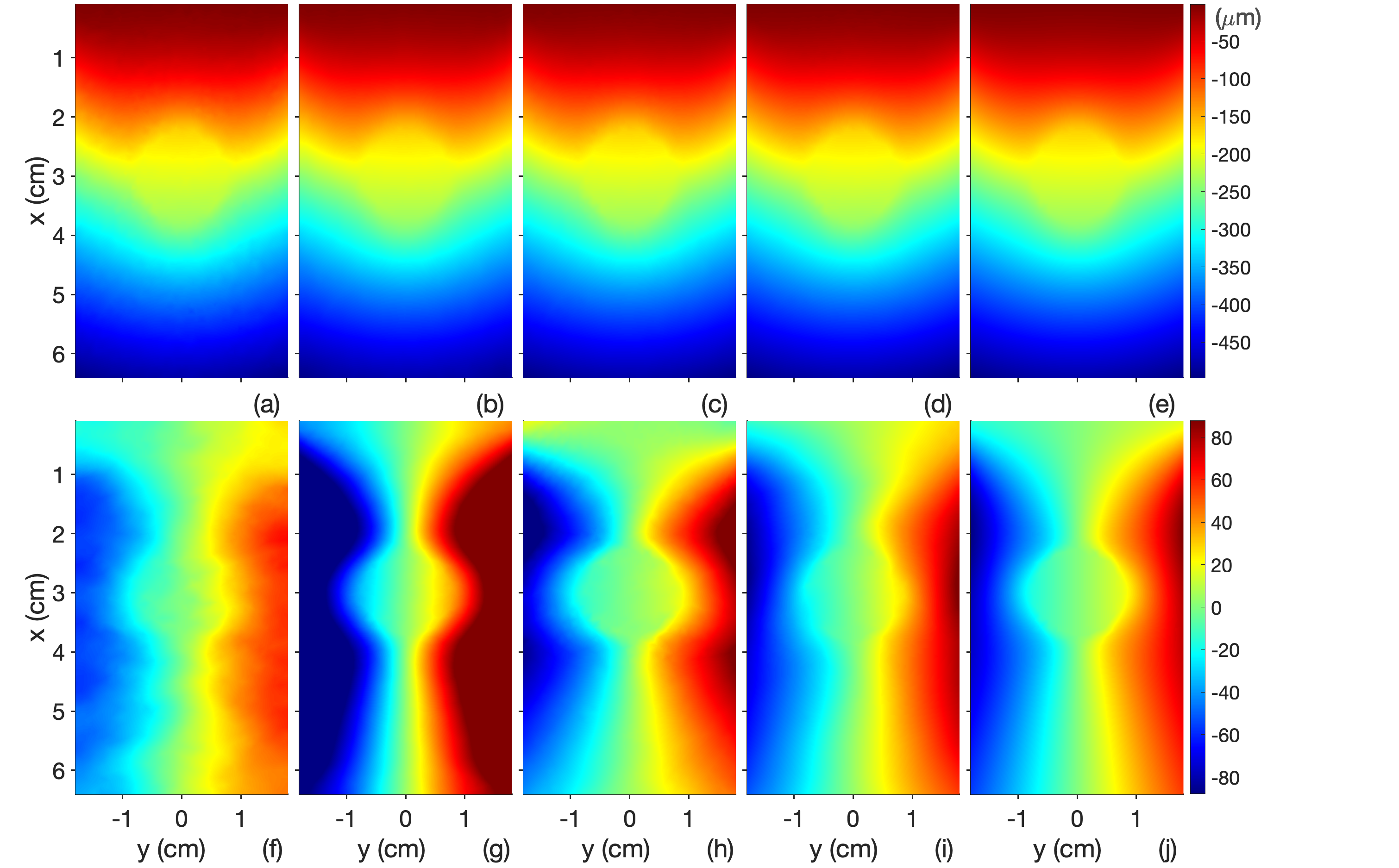}};
        \node[above, left] at (-5.1, 5.25) {$\mathcal{R}_\epsilon$};
	\node[above, left] at (-2.4, 5.25) {$\mathcal{R}_{\epsilon i}$};
	\node[above, left] at (0.3, 5.25) {$\mathcal{R}_{P\epsilon}$};
	\node[above, left] at (2.9, 5.25) {$\mathcal{R}_{P\sigma}$};
	\node[above, left] at (5.6, 5.25) {True};	
	\node[above, rotate=90] at (-7.3, 2.9) {Axial $u_x$};
	\node[above, rotate=90] at (-7.3, -2.1) {Lateral $u_y$};
    \end{tikzpicture}
    \caption{(a-d) Single-frame measured axial displacement estimates ($u_x$) from two US images (frames 1 and 7) created using the 3D forward model (HI) and regularization (a) $\R_\epsilon$, (b) $\R_{\epsilon i}$, (c) $\R_{P\epsilon}$, and (d) $\R_{P\sigma}$. (e) True axial displacement interpolated at nodal location directly from the outputs of the 3D forward model (HI). (f-i) Measured lateral displacement estimates ($u_y$) from 3D forward model (HI) and regularization (f) $\R_\epsilon$, (g) $\R_{\epsilon i}$, (h) $\R_{P\epsilon}$, and (i) $\R_{P\sigma}$. (j) True lateral displacement for the 3D forward model (HI).}
    \label{fig1}
\end{figure}
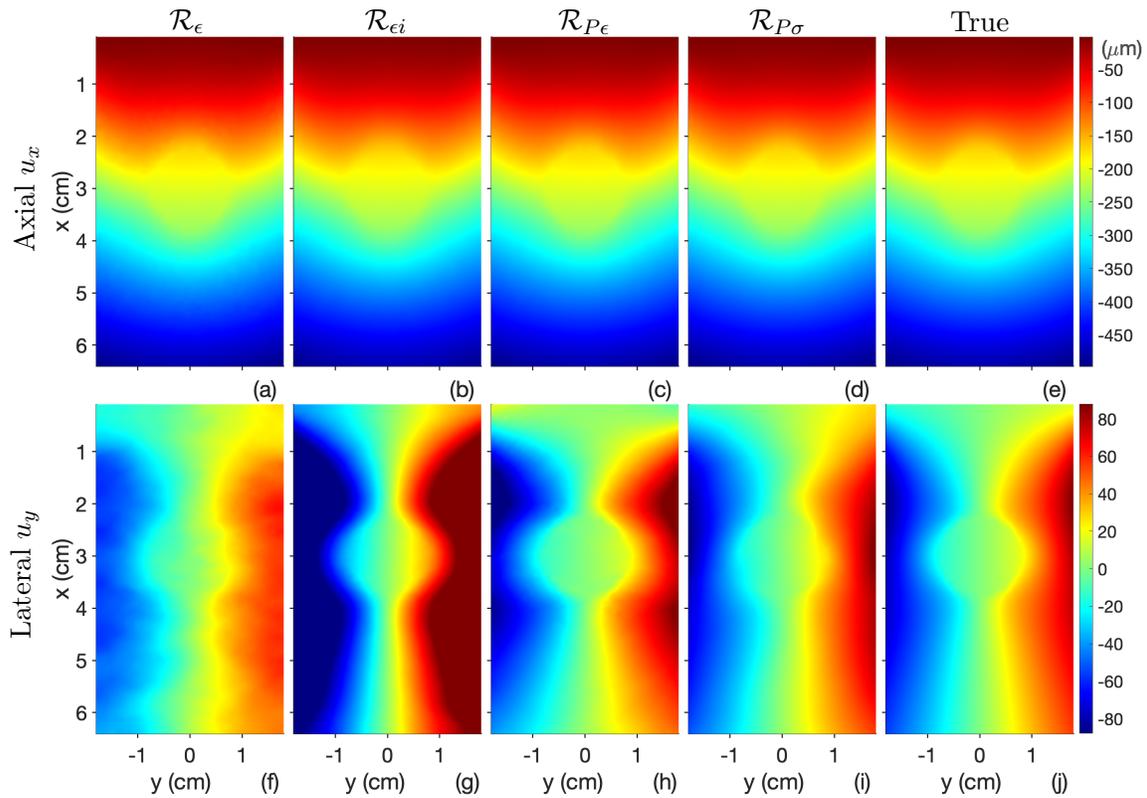
\fi
\ifnum \Plots>0
\begin{figure}[h]
    \begin{tikzpicture}
        \node (image) at (0,0) {\includegraphics[trim={0.1cm 0 2cm 0},clip,width = \textwidth]{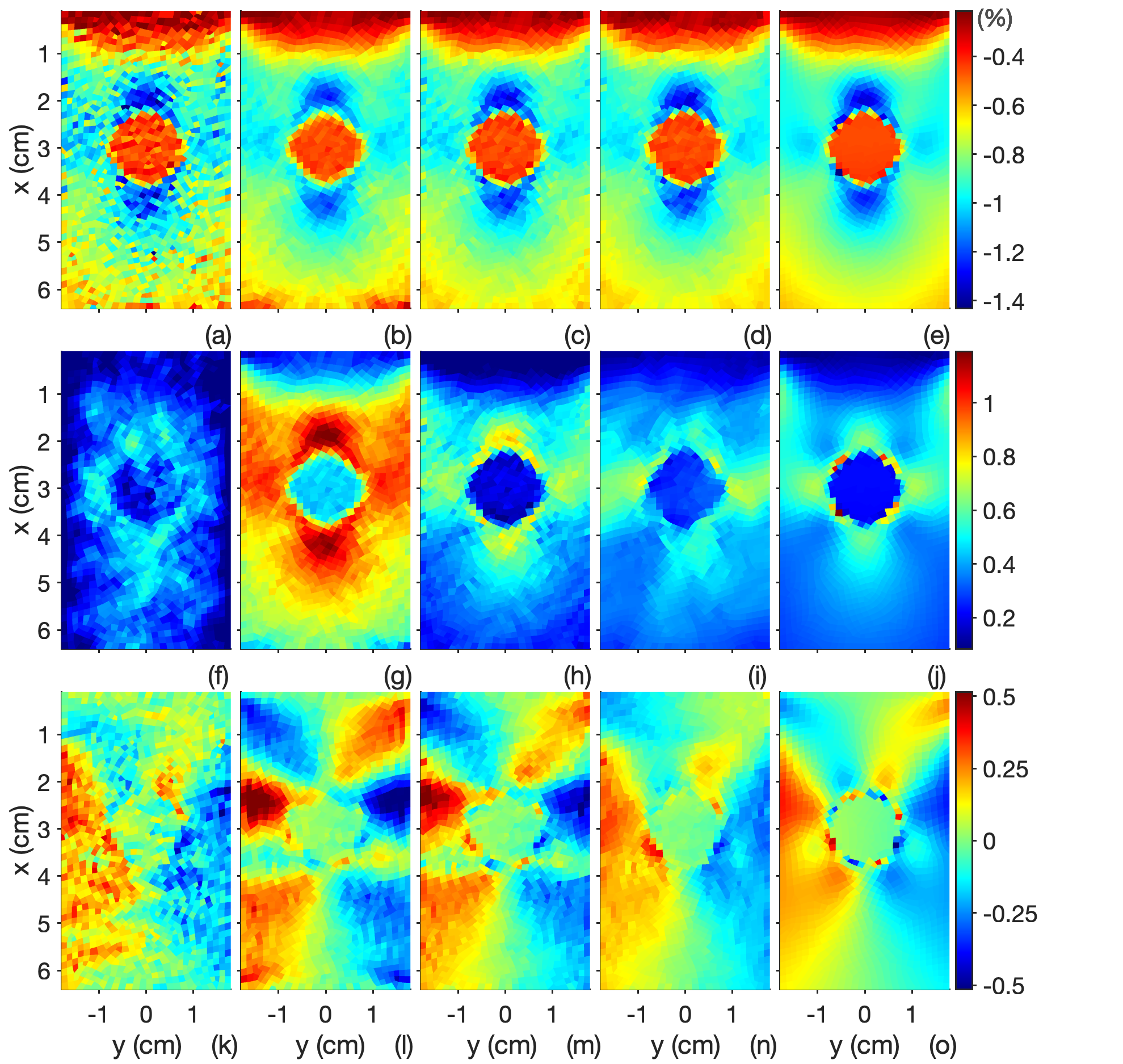}};
        \node[above, left] at (-5.1, 7.75) {$\mathcal{R}_\epsilon$};
	\node[above, left] at (-2.5, 7.75) {$\mathcal{R}_{\epsilon i}$};
	\node[above, left] at (0.2, 7.75) {$\mathcal{R}_{P\epsilon}$};
	\node[above, left] at (2.7, 7.75) {$\mathcal{R}_{P\sigma}$};
	\node[above, left] at (5.2, 7.75) {True};	
	\node[above, rotate=90] at (-7.5, 5.5) {Axial $\epsilon_{xx}$};
	\node[above, rotate=90] at (-7.5, 0.5) {Lateral $\epsilon_{yy}$};
	\node[above, rotate=90] at (-7.5, -4.4) {Shear $\epsilon_{xy}$};
    \end{tikzpicture}
    \caption{(a-d) Measured axial strain ($\epsilon_{xx}$) from two US images (frames 1 and 7) created using the 3D forward model (HI) and regularization (a) $\R_\epsilon$, (b) $\R_{\epsilon i}$, (c) $\R_{P\epsilon}$, and (d) $\R_{P\sigma}$. (e) True axial strain for the 3D forward model (HI). (f-i) Single-frame measured lateral strain ($\epsilon_{yy}$) from two US images using the 3D forward model and regularization (f) $\R_\epsilon$, (g) $\R_{\epsilon i}$, (h) $\R_{P\epsilon}$, and (i) $\R_{Ps\sigma}$. (j) True lateral strain for the 3D forward model (HI). (k-n) Single-frame measured shear strain ($\epsilon_{xy}$) from two US images using the 3D forward model and regularization (k) $\R_\epsilon$, (l) $\R_{\epsilon i}$, (m) $\R_{P\epsilon}$, and (n) $\R_{P\sigma}$. (o) True shear strain for the 3D forward model (HI).}
    \label{fig2}
\end{figure}
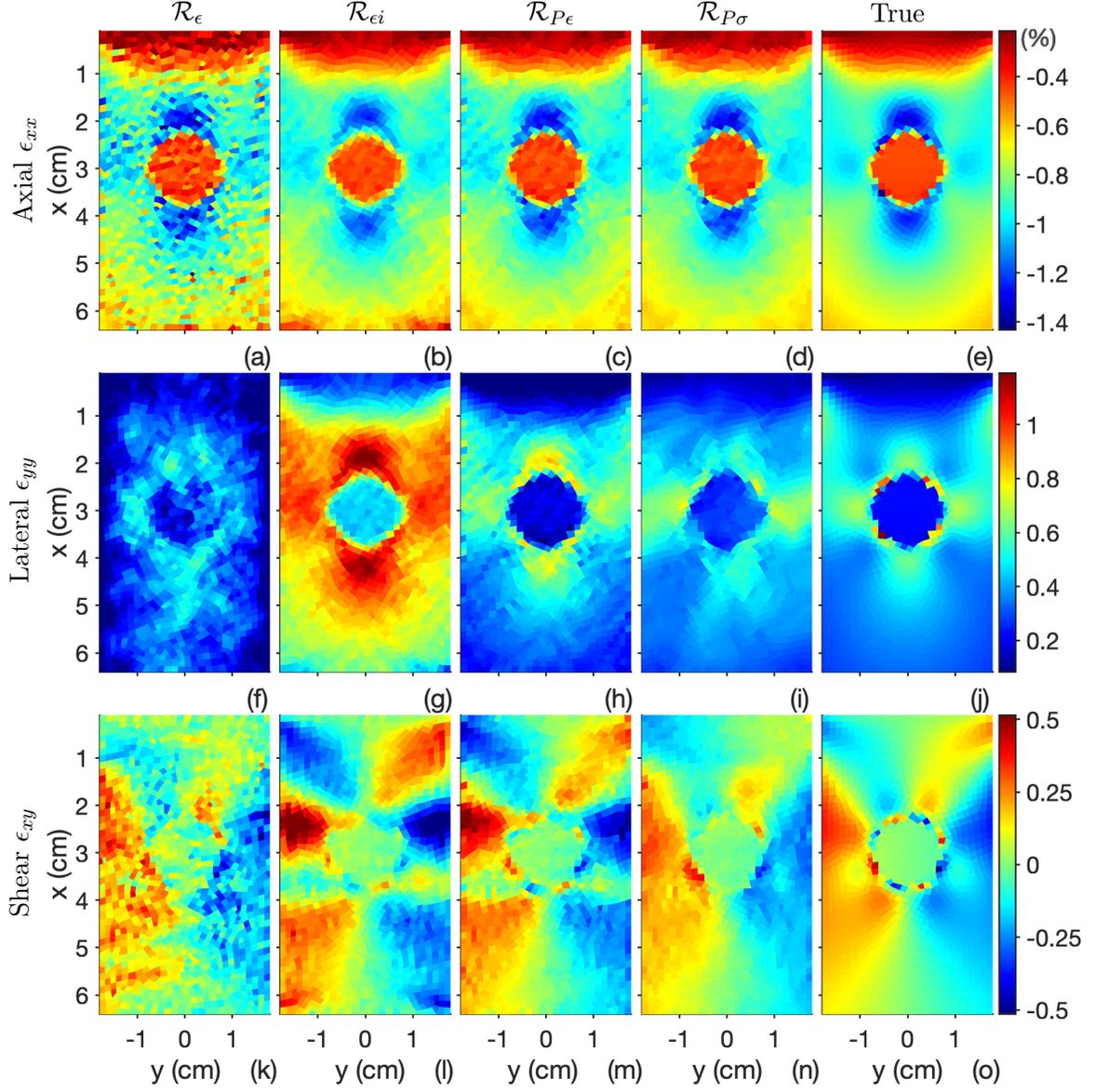
\fi
\ifnum \Plots>0
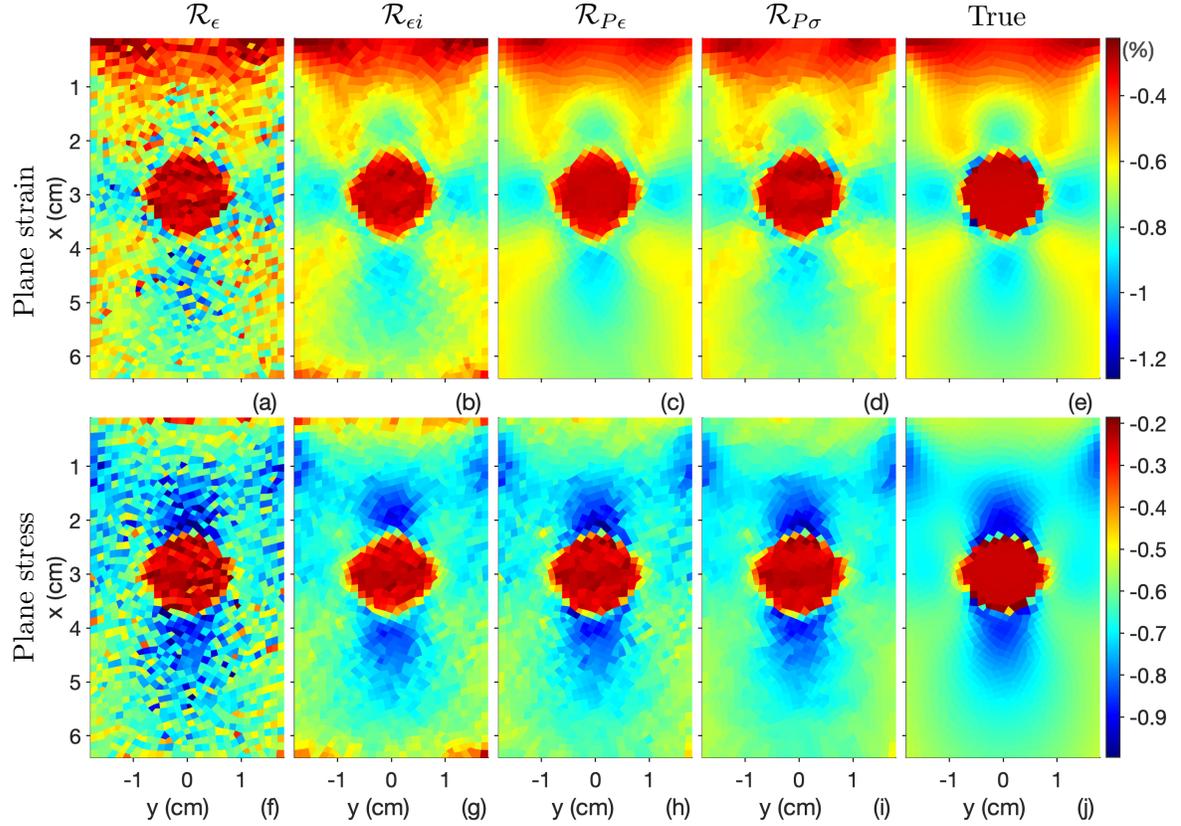
\begin{figure}[h]
    \begin{tikzpicture}
        \node (image) at (0,0) {\includegraphics[trim={1cm 0 4.5cm 0},clip,width = \textwidth]{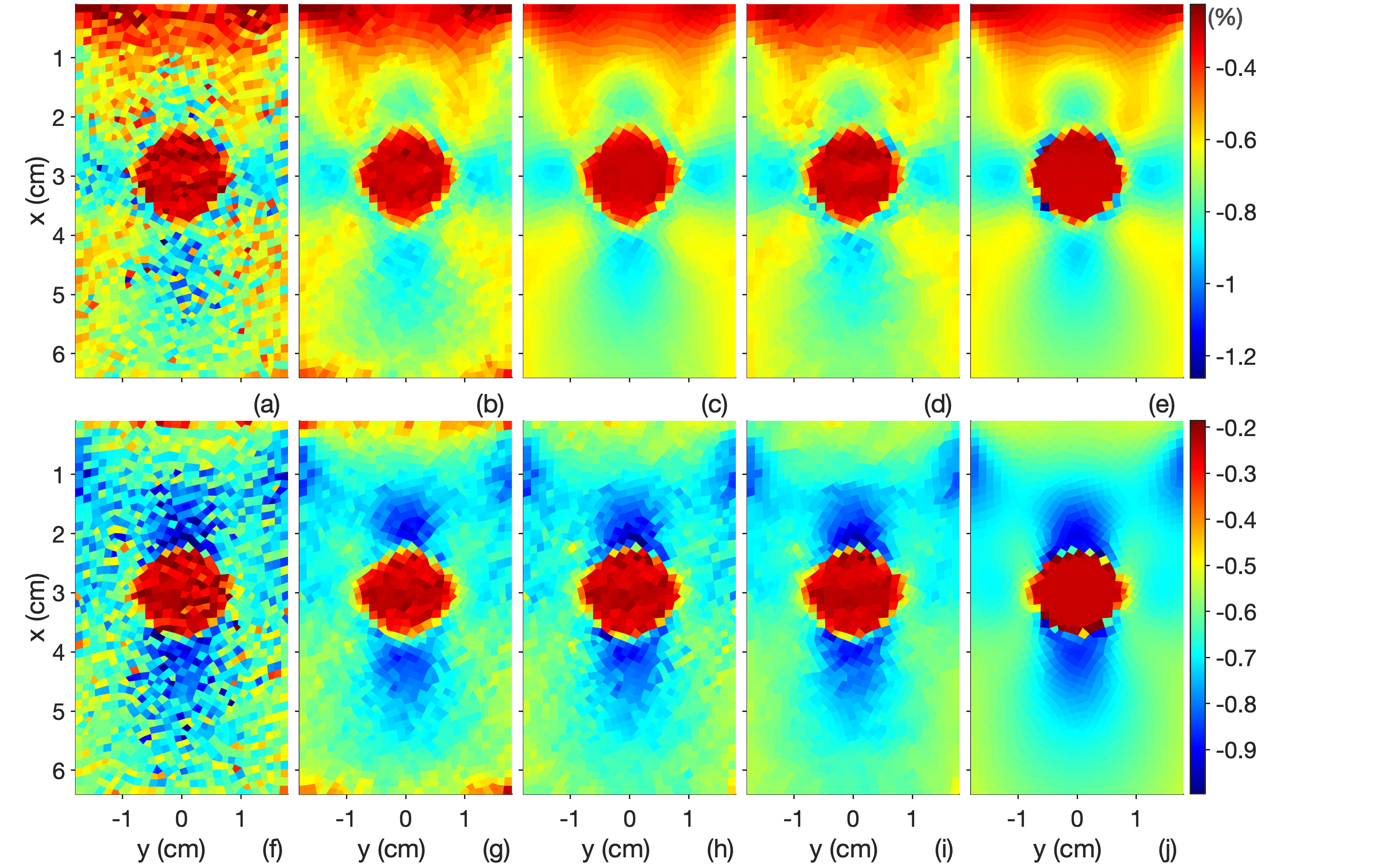}};
        \node[above, left] at (-5.1, 5.5) {$\mathcal{R}_\epsilon$};
	\node[above, left] at (-2.4, 5.5) {$\mathcal{R}_{\epsilon i}$};
	\node[above, left] at (0.3, 5.5) {$\mathcal{R}_{P\epsilon}$};
	\node[above, left] at (2.9, 5.5) {$\mathcal{R}_{P\sigma}$};
	\node[above, left] at (5.6, 5.5) {True};	
	\node[above, rotate=90] at (-7.6, 2.5) {Plane strain};
	\node[above, rotate=90] at (-7.6, -2.1) {Plane stress};	
    \end{tikzpicture}
    \caption{(a-d) Single-frame measurement axial strain ($\epsilon_{xx}$) from two US images (frame 1 and 7) created using the 2D, plane strain forward model (HI) and regularization (a) $\R_\epsilon$, (b) $\R_{\epsilon i}$, (c) $\R_{Ps\epsilon}$, and (d) $\R_{Ps\sigma}$. (e) True axial strain image for the 2D, plane strain forward model (HI). (f-i) Single-frame measured axial strain ($\epsilon_{xx}$) from two US images (frame 1 and 7) created using the 2D, plane stress forward model (HI) and regularization (a) $\R_\epsilon$, (b) $\R_{\epsilon i}$, (c) $\R_{P\epsilon}$, and (d) $\R_{P\sigma}$. (e) True axial strain for the 2D, plane stress forward model (HI).}
    \label{fig3}
\end{figure}
\fi
\ifnum \Plots>0
\begin{figure}[h]
\begin{tikzpicture}
        \node (image) at (0,0) {\includegraphics[trim={1cm 0 4.5cm 0},clip,width = \textwidth]{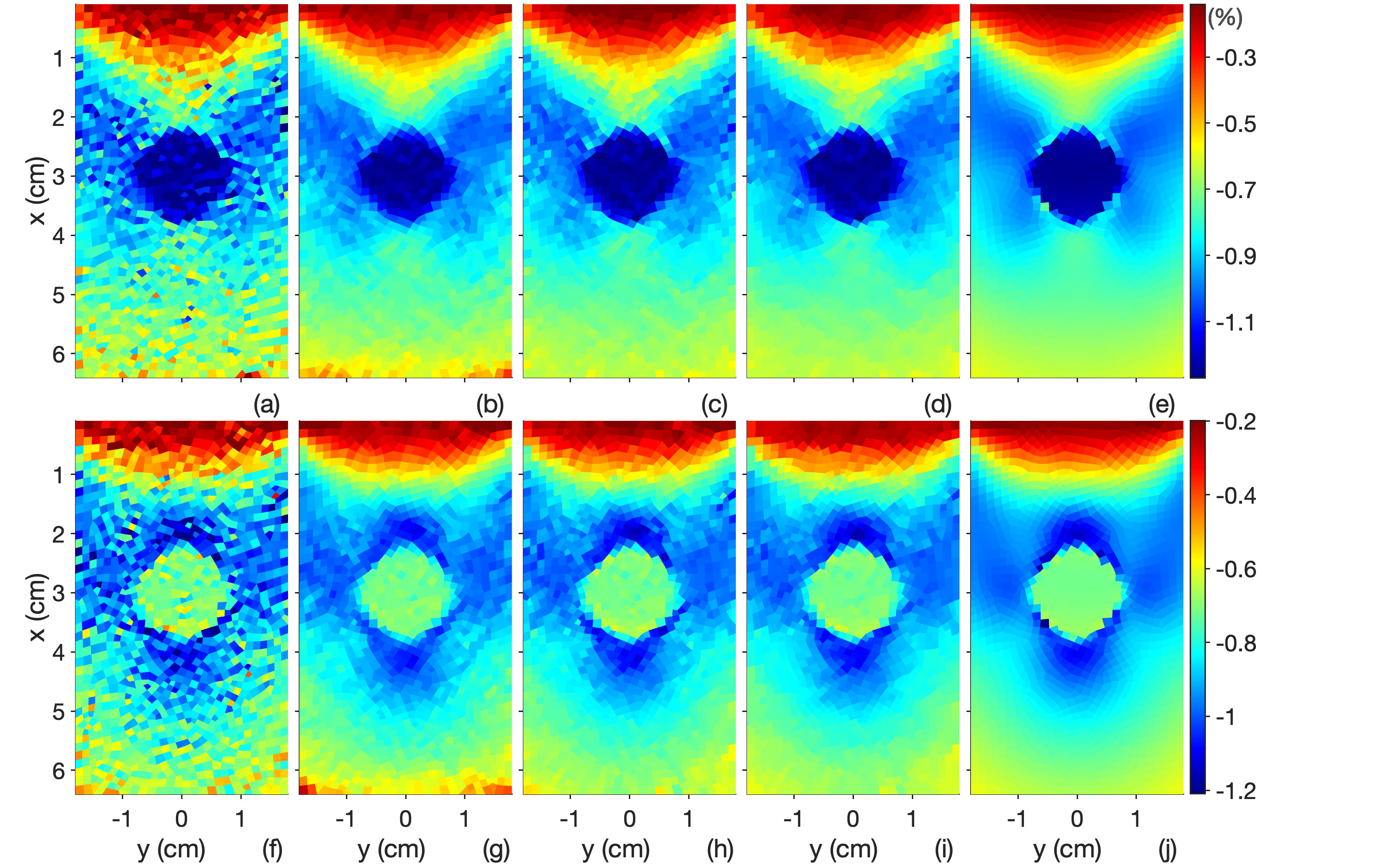}};
        \node[above, left] at (-5.1, 5.5) {$\mathcal{R}_\epsilon$};
	\node[above, left] at (-2.4, 5.5) {$\mathcal{R}_{\epsilon i}$};
	\node[above, left] at (0.3, 5.5) {$\mathcal{R}_{P\epsilon}$};
	\node[above, left] at (2.9, 5.5) {$\mathcal{R}_{P\sigma}$};
	\node[above, left] at (5.6, 5.5) {True};	
	\node[above, rotate=90] at (-7.55, 3) {Soft Inclusion (SI)};
	\node[above, rotate=90] at (-7.55, -2.1) {Med. Hard Inclusion (mHI)};	
    \end{tikzpicture}
    \caption{(a-d) Single-frame measured axial strain ($\epsilon_{xx}$) from two US images (frame 1 and 7) created using the 3D forward model (SI) and regularized as follows (a) $\R_\epsilon$, (b) $\R_{\epsilon i}$, (c) $\R_{P\epsilon}$, and (d) $\R_{P\sigma}$. (e) Exact axial strain image for the 3D forward model (SI). (f-i) Single-frame measured axial strain ($\epsilon_{xx}$) from two US images (frame 1 and 7) created using the 3D forward model (mHI) and regularization (f) $\R_\epsilon$, (g) $\R_{\epsilon i}$, (h) $\R_{P\epsilon}$, and (i) $\R_{P\sigma}$. (j) True axial strain image for the 3D forward model (mHI).}
    \label{fig4}
\end{figure}
\fi

\subsection{Accumulated Displacement Measurements}
The following results correspond to accumulated displacement measurements from frame 1 to frame 20 of the simulated US images. Table \ref{table:acc_strs} shows the calculated values (final frame, $20^{th}$) of the error and contrast metrics for the accumulated displacement measurements of the 3D simulated model sequences with the hardest inclusion (HI). Figures \ref{fig5} and \ref{fig6} show plots of the total displacement and total strain error, respectively, in the accumulated measurements of the 3D simulated model sequences for all regularization types (sub-figures a-d) and inclusion values (colors). For comparison, dual frame measurements are also plotted for measurements of frames 1 to 3 (circles), frames 1 to 7 (squares), and frames 1 to 10 (triangles). Axial and lateral displacement images of the accumulated measurements of all regularization types for the 3D simulated model sequences with the hardest inclusion (HI) are shown in Figure \ref{fig7}. The corresponding axial, lateral and shear strain images are shown in Figure \ref{fig8}.
\begin{table}[h]
\begin{center}
\begin{tabular}{|R{0.65cm}||C{\ctw}|C{\ctw}|C{\ctw}|C{\ctw}|C{\ctw}|C{\ctw}|C{\ctw}|C{0.9cm}|C{1.15cm}|} 
 \hline
 \multicolumn{8}{|c|}{\bf{Measurement Error (Accumulated Disp. Data, HI)}}&\multicolumn{2}{|c|}{\bf{Contrast}}\\
 \hline
 \multicolumn{1}{|c||}{Model/Reg.}&\multicolumn{2}{|c|}{Disp. Comp.} &Total Disp.&\multicolumn{3}{|c|}{Strain Component} &Total Strain&\multicolumn{2}{|c|}{$\epsilon_{xx}$ Only}\\
 \hline\hline
 \multicolumn{1}{|l||}{\bf{3D}}  & $u_{x}$ & $u_{y}$  & $\mathbf{u}$ & $\epsilon_{xx}$ & $\epsilon_{yy}$ & $\epsilon_{xy}$  & $\epsilon$&$SR$&$CNR_e$\\
 \hline\hline
$R_\epsilon$ & $0.1\%$ & $21.4\%$ & $3.1\%$ & $7.5\%$ & $42.8\%$ & $36.0\%$ & $21.3\%$ & $1.68$ & $3.54$ \\
\hline
$R_{\epsilon i}$ & $4.7\%$ & $79.4\%$ & $12.4\%$ & $20.6\%$ & $75.2\%$ & $61.0\%$ & $39.8\%$ & $1.01$ & $6.58e\text{-}4$ \\
\hline
$R_{P\epsilon}$ & $0.1\%$ & $9.6\%$ & $1.4\%$ & $5.1\%$ & $15.7\%$ & $41.1\%$ & $10.8\%$ & $1.68$ & $3.75$ \\
\hline
$R_{P\sigma}$ & $0.1\%$ & $4.4\%$ & $0.6\%$ & $4.0\%$ & $14.4\%$ & $26.1\%$ & $8.6\%$ & $1.68$ & $3.84$ \\
 \hline
\end{tabular}
\caption{Percent displacement and strain error calculated from US images created using the 3D forward model (HI) accumulated over 20 US frames and shown for all regularization types. Results are given for individual strain tensor components and total strain. Contrast recovery is also reported as $SR$ and $CNR_e$. The strain ratio calculated from the true strain of the 3D forward model is $SR = 1.65$ and the contrast to noise ratio is $CNR_e = 3.77$. }
\label{table:acc_strs}
\end{center}
\end{table}

\ifnum \Plots>0
\begin{figure}[h]
\begin{tikzpicture}
        \node (image) at (0,0) {\includegraphics[trim={0cm 0 0cm 0},clip,width = \textwidth]{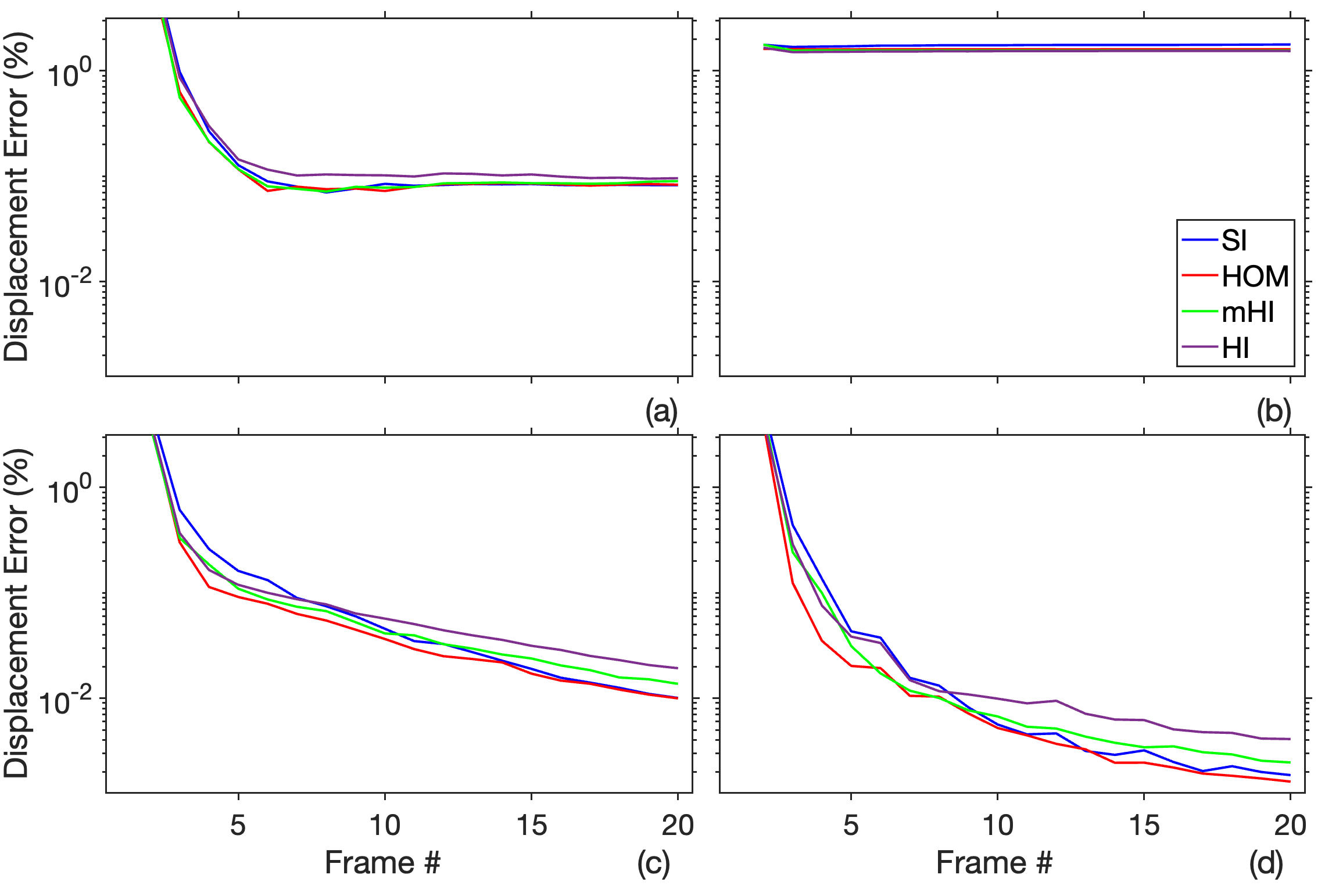}};
        \node[above, left] at (-5.1, 1.5) {$\mathcal{R}_\epsilon$};
	\node[above, left] at (2, 1.5) {$\mathcal{R}_{\epsilon i}$};
	\node[above, left] at (-5.1, -3.2) {$\mathcal{R}_{P\epsilon}$};
	\node[above, left] at (2, -3.2) {$\mathcal{R}_{P\sigma}$};
    \end{tikzpicture}
    \caption{(a-d) Plots of percent error in total displacement ($\mathbf{u}$) estimates from US images created using the 3D forward model with accumulation over 20 US frames with regularization (a) $\R_\epsilon$, (b) $\R_{\epsilon i}$, (c) $\R_{P\epsilon}$, and (d) $\R_{P\sigma}$. Colors correspond to inclusion modulus values of E = 40 MPa (HI, blue), E = 20 MPa (mHI, red), E = 10 MPa (HOM, green), and E = 5 MPa (SI, purple).}
    \label{fig5}
\end{figure}
\fi

\ifnum \Plots>0
\begin{figure}[h]
\begin{tikzpicture}
        \node (image) at (0,0) {\includegraphics[trim={0cm 0 0cm 0},clip,width = \textwidth]{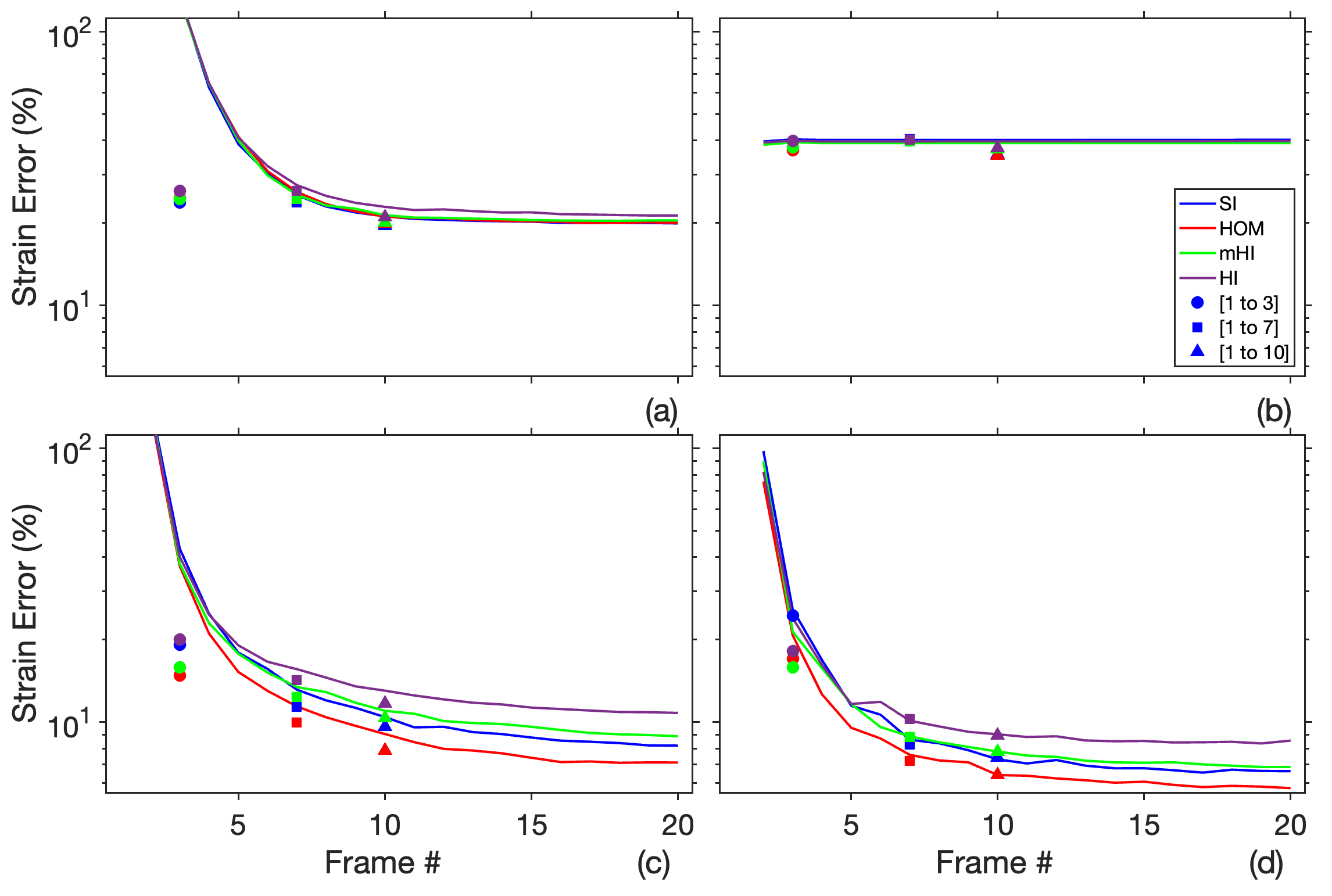}};
        \node[above, left] at (-5.1, 1.5) {$\mathcal{R}_\epsilon$};
	\node[above, left] at (2, 1.5) {$\mathcal{R}_{\epsilon i}$};
	\node[above, left] at (-5.1, -3.2) {$\mathcal{R}_{P\epsilon}$};
	\node[above, left] at (2, -3.2) {$\mathcal{R}_{P\sigma}$};
    \end{tikzpicture}
    \caption{(a-d) Plots of percent error in total strain ($\epsilon$) estimates from US images created using the 3D forward model with accumulation over 20 US frames with regularization (a) $\R_\epsilon$, (b) $\R_{\epsilon i}$, (c) $\R_{P\epsilon}$, and (d) $\R_{P\sigma}$. Colors correspond to inclusion modulus values of E = 40 MPa (HI, blue), E = 20 MPa (mHI, red), E = 10 MPa (HOM, green), and E = 5 MPa (SI, purple). Markers represent error single-frame registration measurements (i.e., 2 US frames,  not accumulated) for frames 1 to 3 (dots), 1 to 7 (squares) and 1 to 10 (triangles).}
    \label{fig6}
\end{figure}
\fi

\ifnum \Plots>0
\begin{figure}[h]
    \begin{tikzpicture}
        \node (image) at (0,0) {\includegraphics[trim={1cm 0 5cm 0},clip,width = \textwidth]{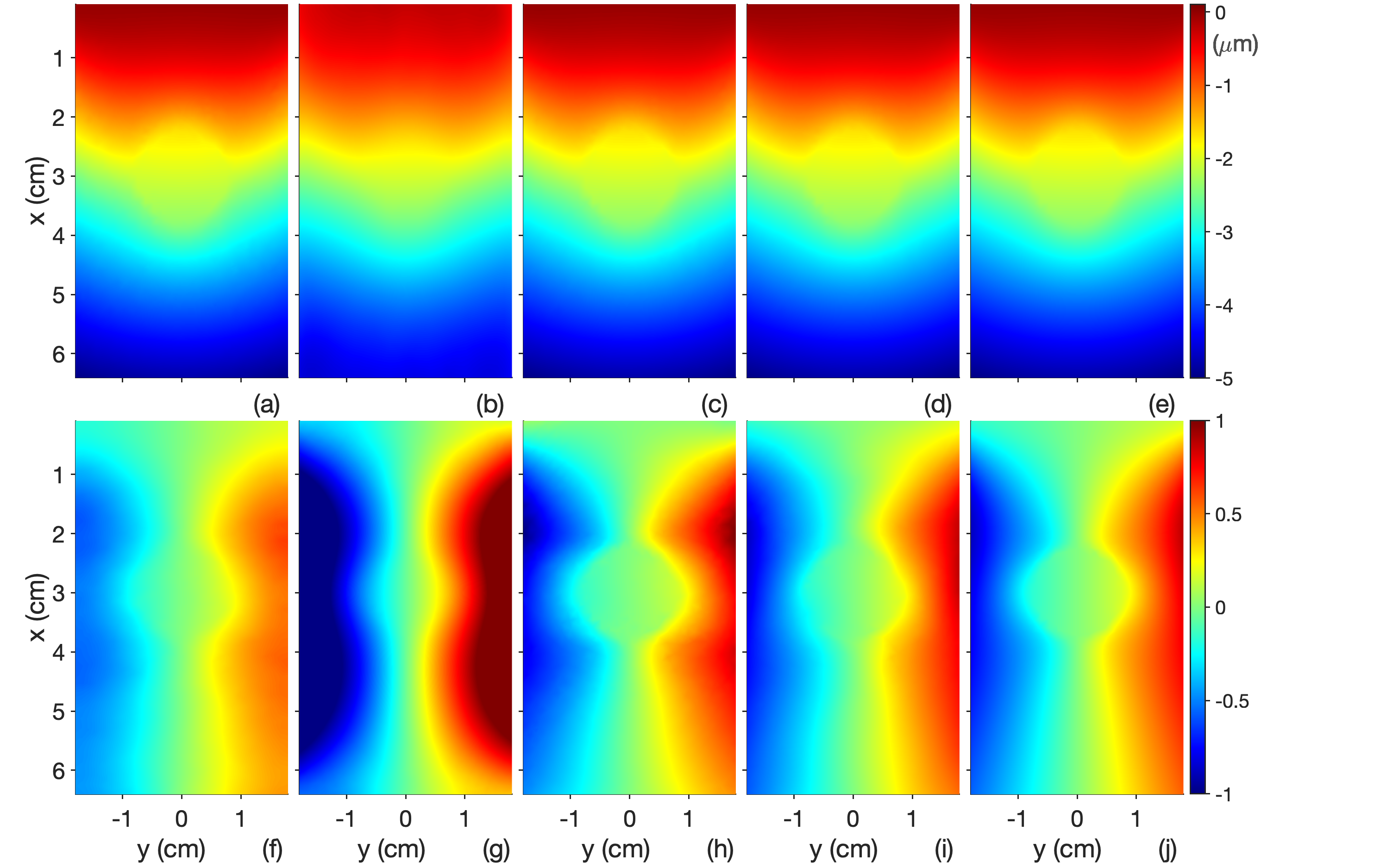}};
        \node[above, left] at (-5.1, 5.5) {$\mathcal{R}_\epsilon$};
	\node[above, left] at (-2.4, 5.5) {$\mathcal{R}_{\epsilon i}$};
	\node[above, left] at (0.3, 5.5) {$\mathcal{R}_{P\epsilon}$};
	\node[above, left] at (3, 5.5) {$\mathcal{R}_{P\sigma}$};
	\node[above, left] at (5.8, 5.5) {True};	
	\node[above, rotate=90] at (-7.6, 3.2) {Axial $u_x$};
	\node[above, rotate=90] at (-7.6, -2.1) {Lateral $u_y$};
    \end{tikzpicture}

    \caption{(a-d) Axial displacement ($u_x$) measurements from US images created using the 3D forward model (HI) with accumulation over 20 US frames for regularizations as follows (a) $\R_\epsilon$, (b) $\R_\epsilon i$, (c) $\R_{P\epsilon}$, and (d) $\R_{P\sigma}$. (e) True axial displacement for the 3D forward model (HI) at 20\textsuperscript{th} frame. (f-j) Lateral displacement ($u_y$) measurements from US images created using the 3D forward model (HI) with accumulation over 20 US frames for regularization (f) $\R_\epsilon$, (g) $\R_{\epsilon i}$, (h) $\R_{P\epsilon}$, and (i) $\R_{P\sigma}$. (j) True lateral displacement for the 3D forward model (HI) at 20\textsuperscript{th} frame.}
    \label{fig7}I
\end{figure}
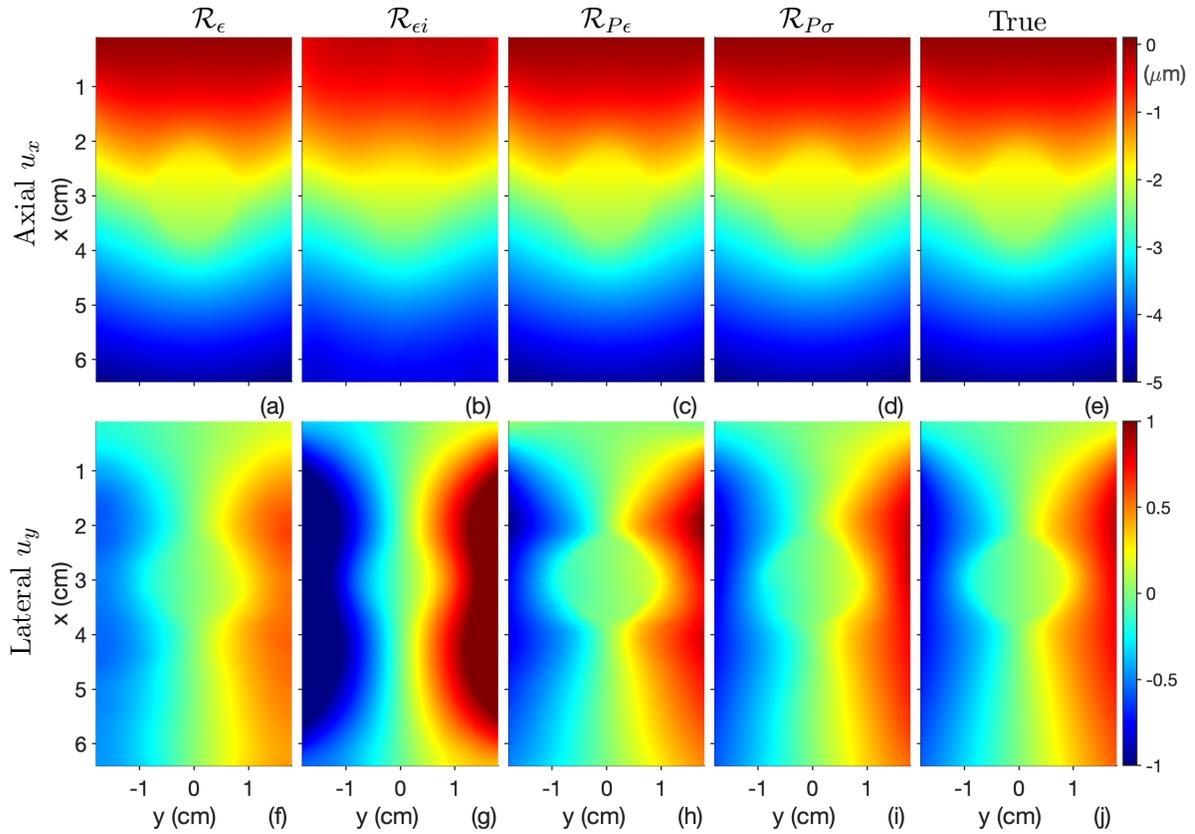
\fi

\ifnum \Plots>0
\begin{figure}[h]
    \begin{tikzpicture}
        \node (image) at (0,0) {\includegraphics[trim={0.25cm 0 3cm 0},clip,width = \textwidth]{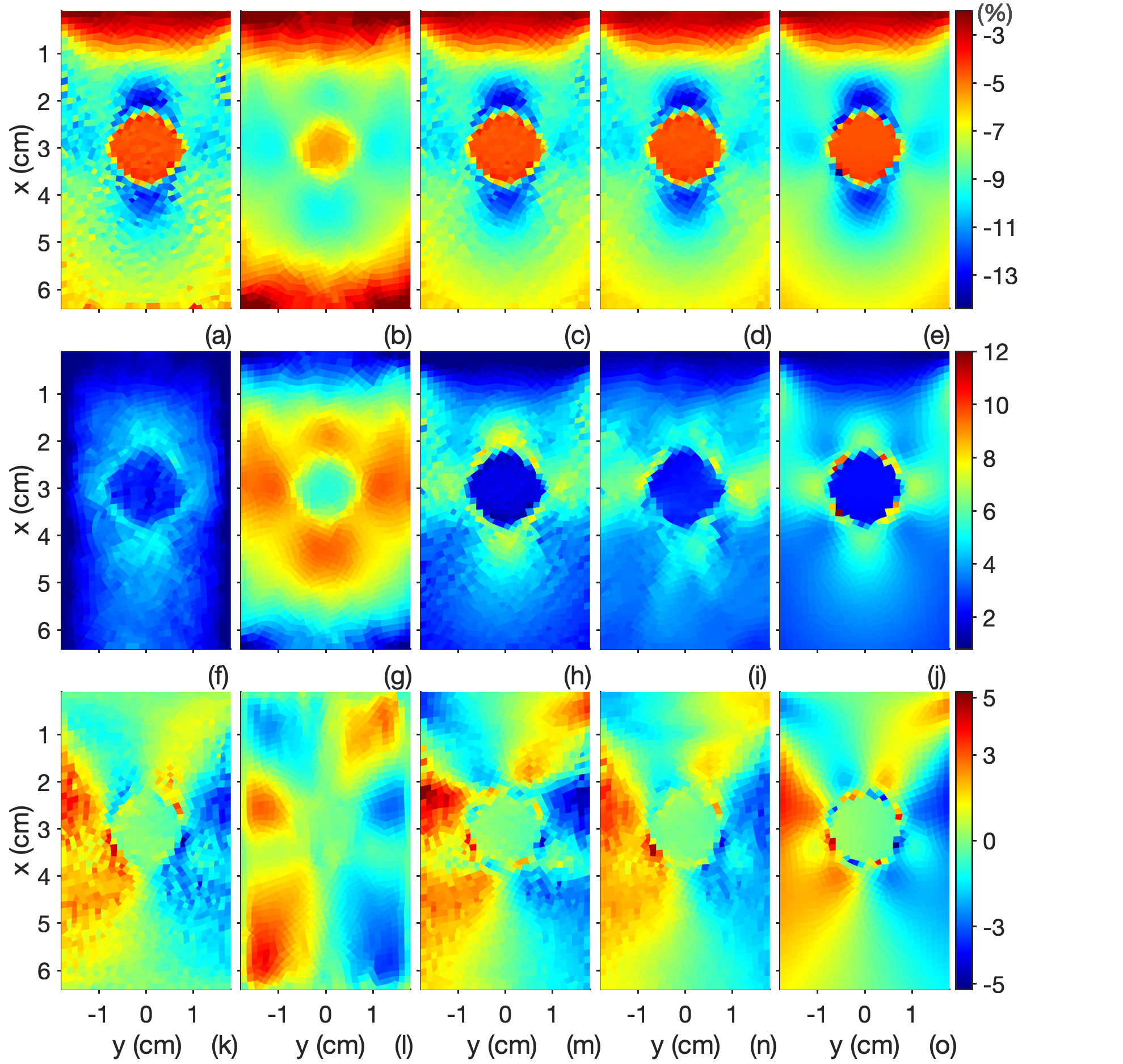}};
        \node[above, left] at (-5.1, 8) {$\mathcal{R}_\epsilon$};
	\node[above, left] at (-2.5, 8) {$\mathcal{R}_{\epsilon i}$};
	\node[above, left] at (0.3, 8) {$\mathcal{R}_{P\epsilon}$};
	\node[above, left] at (3, 8) {$\mathcal{R}_{P\sigma}$};
	\node[above, left] at (5.8, 8) {True};	
	\node[above, rotate=90] at (-7.5, 5.5) {Axial $\epsilon_{xx}$};
	\node[above, rotate=90] at (-7.5, 0.5) {Lateral $\epsilon_{yy}$};
	\node[above, rotate=90] at (-7.5, -4.4) {Shear $\epsilon_{xy}$};
    \end{tikzpicture}        
    \caption{(a-d) Measured axial strain ($\epsilon_{xx}$) from US images created using the 3D forward model (HI) accumulated over 20 US frames and regularized as follows (a) $\R_\epsilon$, (b) $\R_{\epsilon i}$, (c) $\R_{P\epsilon}$, and (d) $\R_{P\sigma}$. (e) True axial strain for the 3D forward model (HI) at 20\textsuperscript{th} frame. (f-i) Measured lateral strain ($\epsilon_{yy}$) from 3D forward model and regularization (f) $\R_\epsilon$, (g) $\R_{\epsilon i}$, (h) $\R_{P\epsilon}$, and (i) $\R_{Ps\sigma}$. (j) True lateral strain for the 3D forward model (HI) at 20\textsuperscript{th} frame. (k-n) Measured shear strain ($\epsilon_{xy}$) from 3D forward model and regularization (k) $\R_\epsilon$, (l) $\R_{\epsilon i}$, (m) $\R_{P\epsilon}$, and (n) $\R_{P\sigma}$. (o) True shear strain for the 3D forward model (HI).}
    \label{fig8}
\end{figure}
\fi

\subsection{Phantom Measurements}
Axial and lateral displacement images of the accumulated measurements of all regularization types for the experimentally collected phantom image sequences are shown in Figure \ref{fig10}. The corresponding axial, lateral and shear strain images are shown in Figure \ref{fig11}. The contrast metric values calculated (final frame, $13^{th}$) from the axial strain fields ($\epsilon^m_{xx}$) for each type of regularization are shown in Table \ref{table:pht} %$.\R_\epsilon$: $SR=2.4$ and $CNR_e=3.7$, $\R_{\epsilon i}$: $SR=1.4$ and $CNR_e=0.45$, $\R_{P\epsilon}$: $SR=2.5$ and $CNR_e=5.7$, and $\R_{P\sigma}$: $SR=2.4$ and $CNR_e=5.6$.

\begin{table}[h]
\begin{center}
\begin{tabular}{ |c|c|c|c|c| } 
 \hline
 Metric & $\R_\epsilon$ & $\R_{\epsilon i}$ & $\R_{P\epsilon}$ & $\R_{P\sigma}$ \\ \hline
 $SR$ & 2.4 & 1.4 & 2.5 & 2.4 \\ \hline
 $CNR_e$ & 3.7 & 0.45 & 5.7 & 5.6 \\ 
 \hline
\end{tabular}
\caption{Strain ratio ($
SR$) and contrast-to-noise $CNR_e$ ration for phantom.}
\label{table:pht}
\end{center}
\end{table}
\ifnum \Plots>0
\begin{figure}[h]
    \centering
        \includegraphics[trim={0.25cm 0 0.5cm 1.2cm},clip,width = 0.6\textwidth]{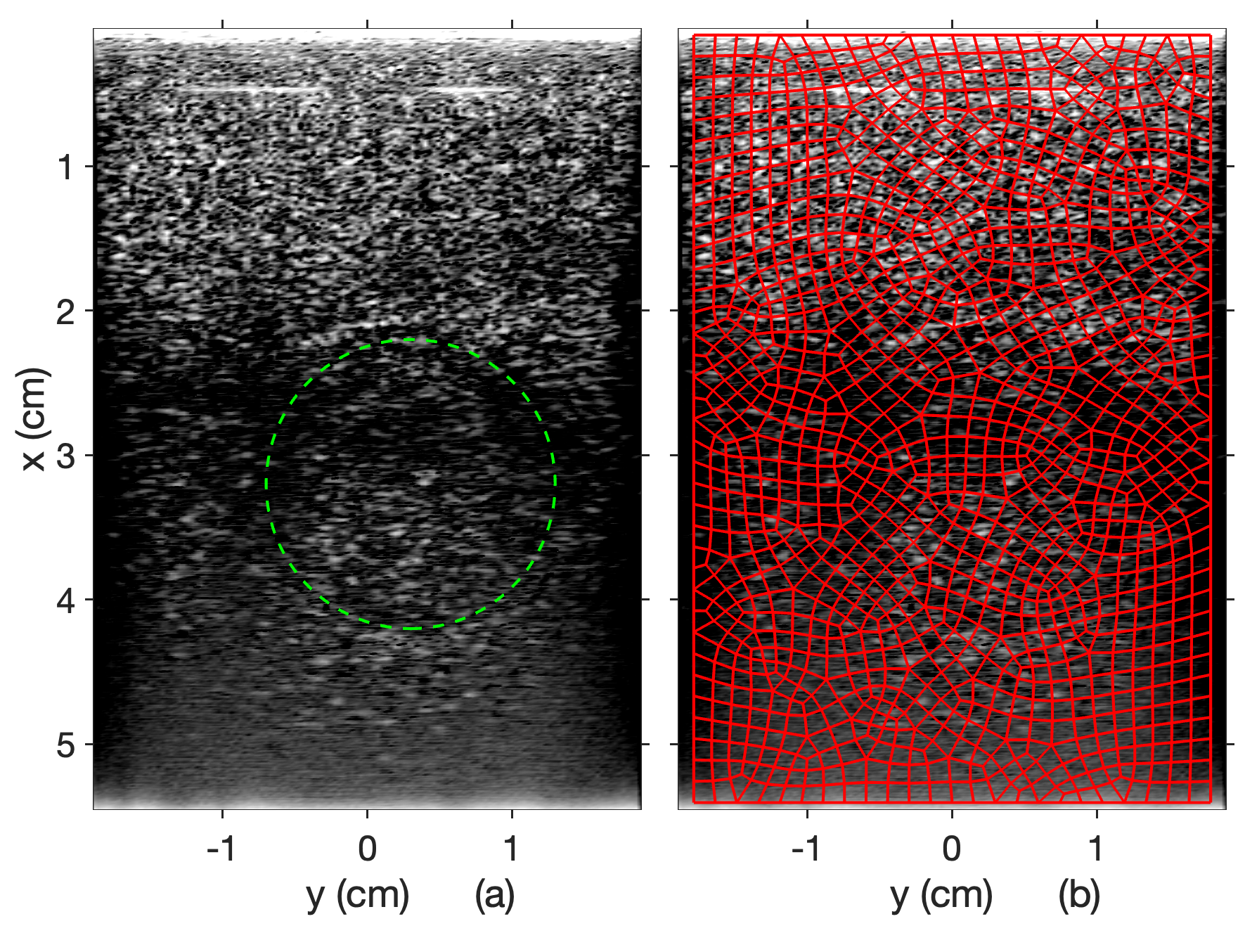}
    \caption{(a) B-mode ultrasound image of elastography tissue mimicking phantom. The green line indicates the region of interest ROI (i.e. the inclusion) used to calculate the strain ratio $SR$ and the contrast-to-noise $CNR_e$ (b) B-mode phantom image with registration mesh (red) overlaid.}
    \label{fig:phantom_og}
\end{figure}
\fi

\ifnum \Plots>0
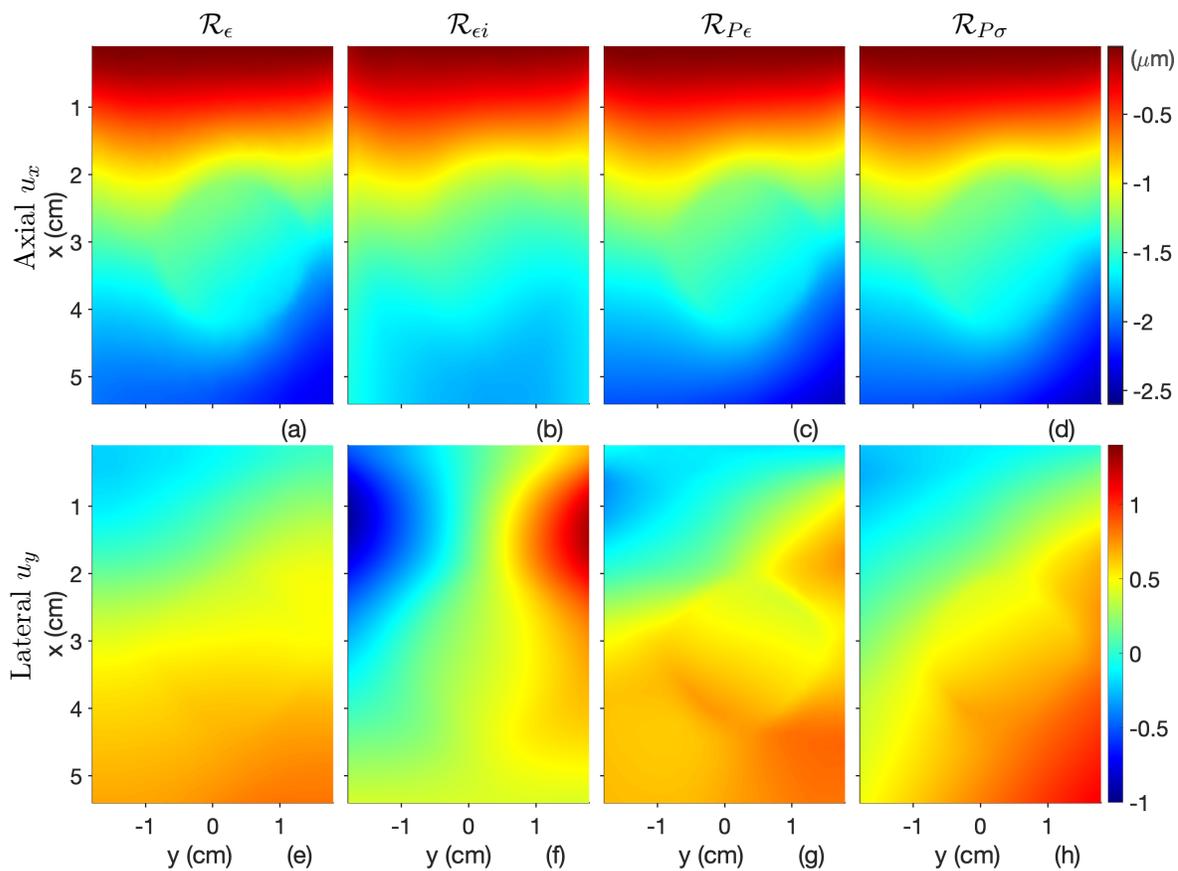
\begin{figure}[h]
    \begin{tikzpicture}
        \node (image) at (0,0) {\includegraphics[trim={0.4cm 0 0.6cm 0},clip,width = \textwidth]{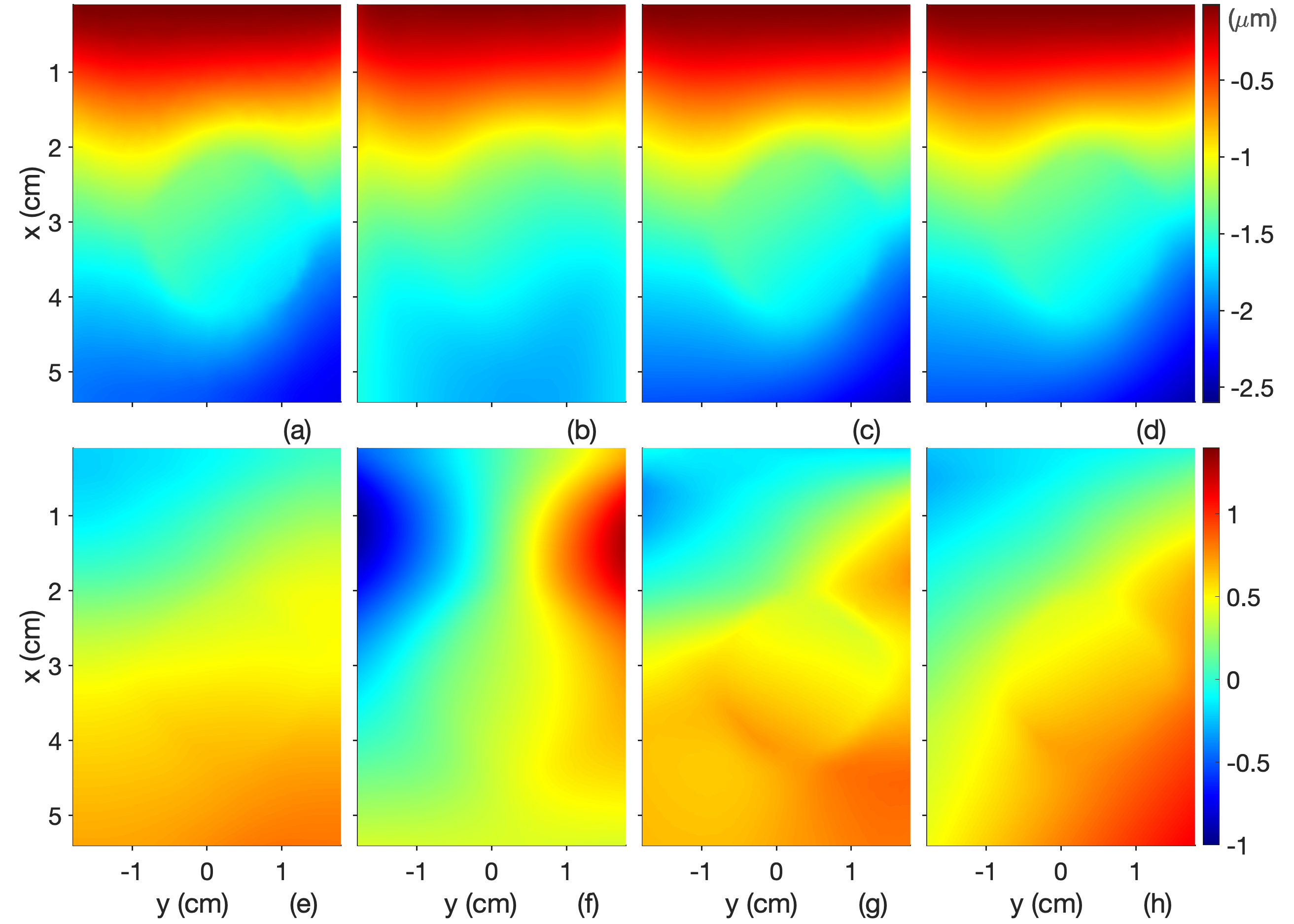}};
        \node[above, left] at (-4.8, 5.75) {$\mathcal{R}_\epsilon$};
	\node[above, left] at (-1.4, 5.75) {$\mathcal{R}_{\epsilon i}$};
	\node[above, left] at (2.1, 5.75) {$\mathcal{R}_{P\epsilon}$};
	\node[above, left] at (5.5, 5.75) {$\mathcal{R}_{P\sigma}$};
	\node[above, rotate=90] at (-7.4, 3.2) {Axial $u_x$};
	\node[above, rotate=90] at (-7.4, -2.1) {Lateral $u_y$};
        \end{tikzpicture}        
    \caption{(a-d) Axial displacement ($u_x$) measurements from US images created using phantom images accumulated over 20 US frames with regularization (a) $\R_\epsilon$, (b) $\R_{\epsilon i}$, (c) $\R_{P\epsilon}$, and (d) $\R_{P\sigma}$. (e-h) Lateral displacement ($u_y$) measurements from US images created using phantom images accumulated over 20 US frames with regularizations as follows (a) $\R_\epsilon$, (b) $\R_\epsilon i$, (c) $\R_{P\epsilon}$, and (d) $\R_{P\sigma}$.}
    \label{fig10}
\end{figure}
\fi

\ifnum \Plots>0
\begin{figure}[h]
    \begin{tikzpicture}
        \node (image) at (0,0) {\includegraphics[trim={0cm 0 0.8cm 0},clip,width = \textwidth]{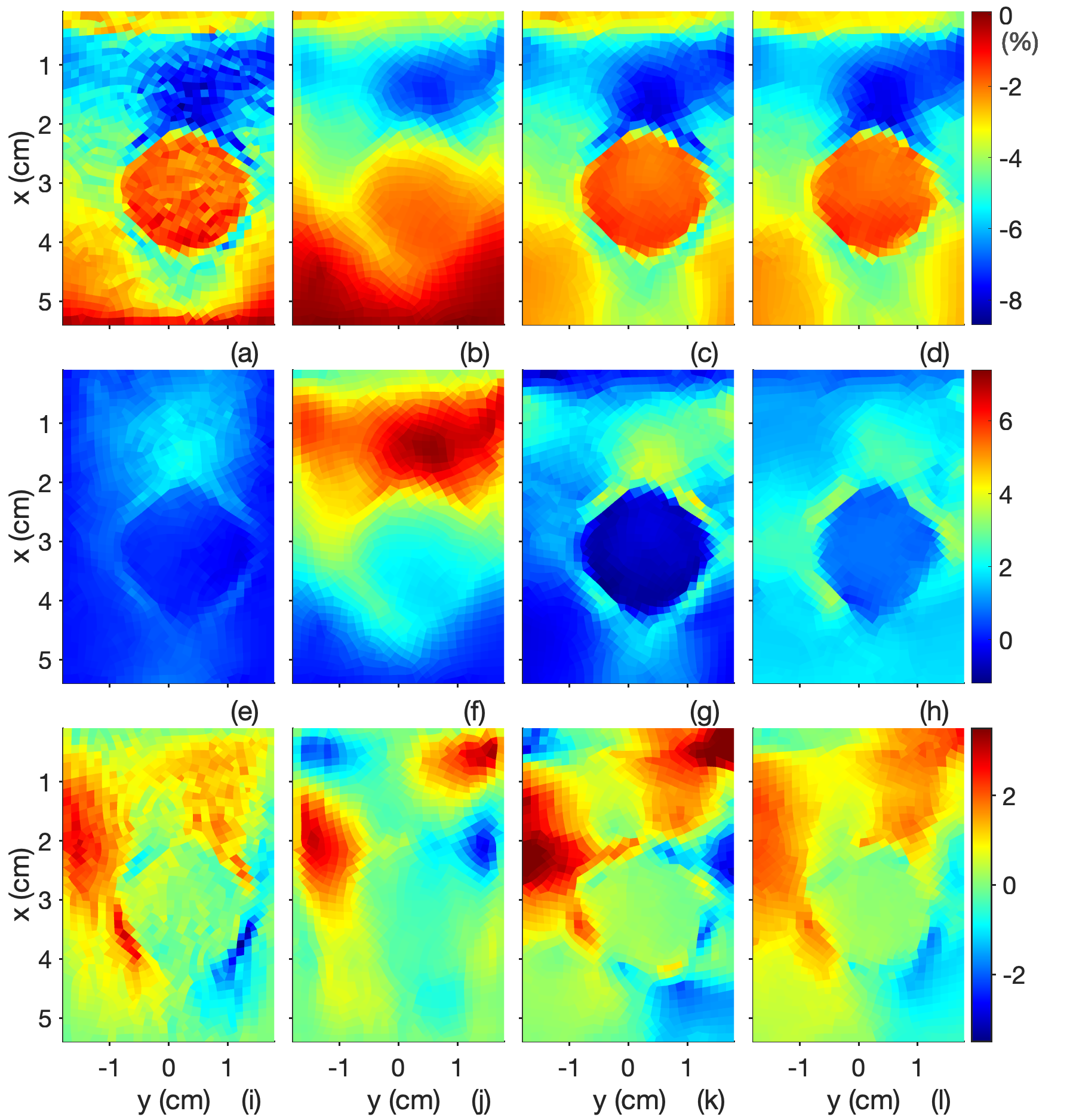}};
        \node[above, left] at (-4.8, 8.25) {$\mathcal{R}_\epsilon$};
	\node[above, left] at (-1.4, 8.25) {$\mathcal{R}_{\epsilon i}$};
	\node[above, left] at (2.1, 8.25) {$\mathcal{R}_{P\epsilon}$};
	\node[above, left] at (5.5, 8.25) {$\mathcal{R}_{P\sigma}$};
	\node[above, rotate=90] at (-7.4, 5.7) {Axial $\epsilon_{xx}$};
	\node[above, rotate=90] at (-7.4, 0.5) {Lateral $\epsilon_{yy}$};
	\node[above, rotate=90] at (-7.4, -4.6) {Shear $\epsilon_{xy}$};
        \end{tikzpicture}        

    \caption{(a-d) Axial strain ($\epsilon_{xx}$) measurements from US images created using phantom images accumulated over 20 US frames with regularization (a) $\R_\epsilon$, (b) $\R_{\epsilon i}$, (c) $\R_{P\epsilon}$, and (d) $\R_{P\sigma}$. (e-h) Lateral strain ($\epsilon_{yy}$) measurements from phantom US images with regularization (a) $\R_\epsilon$, (b) $\R_{\epsilon i}$, (c) $\R_{P\epsilon}$, and (d) $\R_{P\sigma}$. (i-l) Shear Strain ($\epsilon_{xy}$) measurements from phantom US images with regularization (a) $\R_\epsilon$, (b) $\R_{\epsilon i}$, (c) $\R_{P\epsilon}$, and (d) $\R_{P\sigma}$.}
    \label{fig11}
\end{figure}
\fi

\section{Discussion}
This work provides a direct comparison of image registration techniques using a strain-based, strain-based plus incompressibility constraint, and novel mechanics-based regularization strategies for displacement estimation in quasi-static ultrasound elastography. As is the case for most regularization in optimization techniques, the accuracy of the solution can be highly dependent on the choice of the regularization parameter. Since the goal of this work was to compare multiple regularization types, we chose to determine the appropriate regularization parameter for each regularization type by using a single frame-to-frame measurement of two representative simulated images. The images were chosen to differ by approximately the same amount of deformation as the maximum frame-to-frame strain realized in the forward simulation of the multiframe data. Multiple regularization solutions were found for each type of regularization, with varying magnitudes of $\alpha$ and the value that produced the most accurate solution in the strain error was used for the remaining registrations in this work. It is important to note that this method may not be the best or most practical method for choosing the parameter in all situations. Appendix Section \ref{apss:par} provides further details to the process used in this work.

Our results show that, when applied to simulations of two US images of a 3D deforming material, traditional strain regularization ($\R\epsilon$) yields reasonable axial displacement estimates, with average error as low as 0.2$\%$ (see Table \ref{table:strs}). In fact, all regularization types showed comparable levels of error in the axial displacement and strain functions. This may be because, in part, we have very low frame-to-frame strains. However, lateral displacement estimates using traditional strain regularization are as high as $\sim 20\%$ and the results using strain regularization with an incompressibility constraint were as high as $97\%$. The plane strain SPREME regularization improved from strain regularization with $\approx17\%$ displacement error. However, SPREME plane stress regularization outperformed traditional strain regularization by reducing the lateral displacement error by more than half. In addition, the strain-based regularization also underperformed compared to the model-based methods when displacement fields are accumulated over multiple frames.

The error in the strain components of the 3D deforming material showed a similar trend (see Table \ref{table:strs}). Interestingly, the error in the axial and lateral strain estimates was approximately three times less using SPREME plane stress regularization versus traditional strain regularization. In addition, plane stress regularization produced the lowest shear strain errors of $\approx 34\%$, about half that of traditional strain regularization. While the plane strain SPREME regularization also performed better than traditional strain regularization, it did not perform as well as plane stress SPREME regularization in the 3D model. Strain regularization with an incompressibility constraint performed the worst in all strain error metrics compared to the other regularization types. Although this approach reduced the lateral strain error in idealized 2D plane strain simulations, it introduced substantial errors and artifacts in the 3D and accumulated displacement measurements. This reflects a sensitivity to the specific deformation mode and indicates that strong incompressibility enforcement may over-constrain the solution under more realistic conditions.

As expected, the best performance was achieved when the assumptions of the regularization model matched the model used in the forward simulations. Among the techniques evaluated, momentum-based regularization under plane stress assumptions ($\R_{P\sigma}$) consistently produced the most accurate displacement and strain estimates, particularly in the 3D deformation case where this 2D assumption is violated. The simulations where the forward model was indeed plane strain was the only case where the measurement errors of the strain regularization with an incompressibility constraint performed well.

The results of displacement and strain error were largely consistent across all different magnitudes of inclusion tested using the 3D forward model, suggesting little dependence of the regularization outcome on the underlying modulus distribution (see Tables \ref{table:strs} and \ref{table:Astrs}). However, it is true that the improvement in contrast realized in the SPREME plane stress was the most pronounced when the contrast itself was the highest (HI).

For accumulated measurements performed on the 3D image simulations, the plane strain and plane stress model-based regularization schemes were more accurate than the strain with and without incompressibility constraint (see Table \ref{table:acc_strs}). The $SR$ and $CNR_e$ resulting from strain regularization with an incompressibility constraint were significantly lower than those of the other methods, suggesting an accumulation of bias in those measurements. It also appears that the error in the displacement estimates of the SPREME regularization types continues to decrease with increasing frames and strain, where both strain regularization types have plateaued as demonstrated in Figure \ref{fig5}.

Visually, the displacement and strain images of the SPREME regularization types appeared to be less noisy and better match the images of the true simulated displacements and strains in most cases (see Figures \ref{fig1} -- \ref{fig4}, \ref{fig7}, \ref{fig8}). In particular, the SPREME plane stress regularization appeared to be, in all cases, the best qualitatively matching to the true values. This suggests that, when the true underlying deformation is unknown or difficult to control experimentally or clinically, using a SPREME plane stress regularization may be the best option. Visual inspection of the phantom images suggests that the SPREME regularization types also result in the best strain images (see Figure \ref{fig11}). The displacement and strain images of the strain regularization with incompressibility constraint show significantly different results than the other regularization types, further suggesting that an incompressibility constraint introduces a bias that is inconsistent with the phantom generated experimental data (see Figures \ref{fig10}, \ref{fig11}).

Despite the strengths of this comparative analysis, several limitations warrant consideration. First, all regularization strategies were tuned using simulated data with relatively small frame-to-frame strain ($\sim0.8\%$), which may not generalize to higher strain scenarios or pathological tissues with more complex mechanical behavior. Second, the assumption of piecewise constant shear modulus in momentum-based regularization is an idealization that may not hold \textit{in vivo}, particularly near tissue interfaces. Additionally, although the 3D simulations included out-of-plane deformations, the image registration and regularization frameworks were strictly two-dimensional, limiting their ability to capture full volumetric strain patterns. This introduces a mismatch between the model assumptions and the ground truth in the validation data. Lastly, while the phantom data helped demonstrate applicability to experimental conditions, no \textit{in vivo} imaging was performed, and future work is needed to assess robustness in the presence of physiological motion, speckle decorrelation, and real-world image artifacts.

\section{Summary and Conclusions}
This work introduced a SPREME momentum-based regularization that penalizes deformations that are inconsistent with 2D plane strain and plane stress, which is similar to a displacement measurement filter developed by Babaniyi et. al. \cite{babaniyi_spreme}. These techniques were compared directly to strain-based regularization schemes with and without a penalty imposing 2D tissue incompressibility. We showed that the accuracy of SPREME regularization types is higher than strain-based regularization schemes and that imposing 2D tissue incompressibility can introduce significant bias. Our results also show that SPREME regularization is robust to variations in modulus contrast and that measurement accuracy can improve with accumulation. Measurements from real US image data obtained using an US elastography phantom for all regularization types seemed consistent with our findings using simulated data.

The adaptation of SPREME-like momentum regularization into the registration framework ($\R_{P\epsilon}$ and $\R_{P\sigma}$) offers a practical route to enforcing mechanical consistency without requiring knowledge of boundary conditions or material properties. We showed that using our novel SPREME regularization will not only produce better results but will also bias the data toward a model-consistent measurement. Thus, if further processing steps are to be utilized, such as modulus inversion, we would hypothesize that using SPREME regularization with a consistent plane strain or plane stress inversion model would also improve reconstructed modulus images. These measurements are also discretized on an FE mesh as well, avoiding interpolation error if an FE inversion method is utilized. 

Finally, we believe that the results herein demonstrate the advantage of using model-based regularization schemes for US registration applied to static elastography. 

\section*{Acknowledgments}
The authors thank Dr. Paul E. Barbone for his guidance and thoughtful discussions of this work. The authors also acknowledge the use of the tools of ChatGPT 5.0, Grammarly and Writefull to assist in the editing for grammar and wording of this manuscript. This work was supported by the National Science Foundation LEAPS-MPS grant under Grant number 2213493.
\appendix
\section{Appendix}
\subsection{Numerical Integration}
\label{apss:im}
The displacement fields used within the image matching term and all of the regularizations types, are discretized on a mesh consisting of 4 node, bilinear quadrilateral elements. The strain magnitude regularization terms of the strain-based regularizations, $\R_\epsilon[\bu(\bx)]$ and $\R_{\epsilon i}[\bu(\bx)]$, and their derivatives were numerically integrated using a 3x3 Gaussian quadrature method. The incompressibility constraint term in $\R_{\epsilon i}[\bu(\bx)]$ (that is, the second term of Eq.\@ (\ref{eq:saicr})) was integrated with a reduced integration scheme, where integrations were performed using a midpoint method at the centers of finite elements to avoid over constraining the solution \cite{hughes2012finite}. Integration of the momentum-based regularizations terms were treated differently, as explained in  Appendix Section \ref{apss:mbr}.

The image matching functional term of Eq.\@ (\ref{eq:IR}), and it's derivatives, were integrated using a midpoint method at the pixel locations as discrete points. The discretized displacement fields of the FE elements were interpolated to pixel values using an inverse mapping of the bilinear FE shape functions \cite{zhao1999consistent}. The values of the images outside the FE meshed region, $\Omega$, were assumed to be zero.

The interpolation of $I_2(\bx+\bu(\bx))$ in Eq.\@ (\ref{eq:IR}) was performed using the Matlab function ``griddedInterpolant()" with a cubic interpolation. If the magnitude of $\bx+\bu(\bx)$ was large enough, at a given position $\bx$, to fall outside the imaged domain, the value of the interpolation, $I_2$, was set to 0.

\subsection{Reduced Form of TV Momentum Regularization}
\label{apss:mbr}
In the model-based regularization of Eq.\@ (\ref{eq:mreg}), the momentum has an explicit second derivative of $\bu(\bx)$, (i.e., in $\nabla \cdot \bs{A}[\bu(\bx)]$). For FE approximations, these second derivatives cannot be calculated directly. However, if we assume that $\bs{A}[\bu(\bx)]$ is piecewise constant within an element and we use a total variational regularization approach, the derivatives of $\bs{A}[\bu(\bx)]$ can be avoided. A simple example of this reduced form is presented here, where we show that the TV regularization of a linearly varying symmetric tensor function, $\bs{A}[\bu(\bx)]=\bs{A}(\bx)$, will reduce to a sum of the magnitudes of element jumps of a piecewise constant function, in the limit of vanishing elements. Consider a linearly varying $\bs{A}(\bx)$ on a simple, regular $3\times 3$ FE mesh shown in Figure \ref{fig:red}, where $\bs{A}(\bx) = \bs{A}i(\bx)$ is constant within the corner elements $(i = 1 \ \text{to} \ 4)$. Thus, in the center-edge elements, $\bs{A}i(\bx)$ $(i = 5 \ \text{to} \ 8)$ varies linearly between the corner elements. The center element $(i = 9)$, $\bs{A}(\bx)$ is bilinear. When the integral in Eq.\@ (\ref{eq:mer}) is considered for this discretization, it is sufficient to consider the integrals of each type of element individually and a piecewise constant discretization for $\bs{A}(\bx)$ can be recovered by considering the limit of this discretization as the orthogonal sides ($\Delta x$ and $\Delta y$ from Figure \ref{fig:red}) go to zero. 
\ifnum \Plots>0
\begin{figure}[H]
    \centering
    \includegraphics[trim={1.6cm 0.8cm 12cm 0.7cm},clip,width = 8cm]{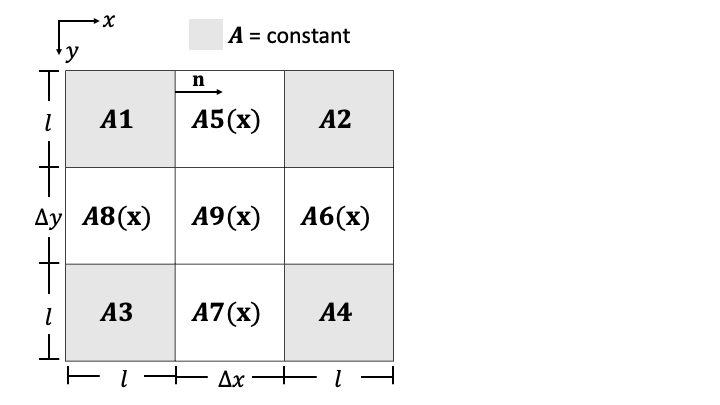}   
    \caption{Schematic demonstrating the simplification of the momentum TV regularization integration.}
    \label{fig:red}
\end{figure}
\fi
When considering the integration of Eq.\@ (\ref{eq:mer}) of this example mesh, the value of $\nabla \cdot \bs{A}(\bx)$ in the corner elements $(i = 1 \ \text{to} \ 4)$ will be exactly 0 and thus can be neglected. For center-edge elements (e.g., $i=5$) the integral in equation (\ref{eq:mreg}) is:
\beq
\int_{\Omega_{5}} ||\nabla\cdot\bs{A}(\bx)||_2 \ d\Omega_{5} = 
\bigintsss_{\Omega_{5}} \sqrt{\left(\frac{\partial A_{xx}}{\partial x}\right)^2+\left(\frac{\partial A_{xy}}{\partial x}\right)^2} \ d\Omega_{5},
\eeq
where partial derivatives of $\bs{A}$ w.r.t. $y$ are $0$ for element $i=5$. Furthermore, the partial derivatives w.r.t $x$ are constant within these elements, thus the integral further reduces to
\beq
\int_{\Omega_{5}} ||\nabla\cdot\bs{A}(\bx)||_2 \ d\Omega_{5} = (l \ \Delta x) \sqrt{\left(\frac{A2_{xx} - A1_{xx}}{\Delta x}\right)^2+\left(\frac{A2_{xy} - A1_{xy}}{\Delta x}\right)^2} 
\label{eq:Ae5}
\eeq
where $(l \ \Delta x)$ is the area of the square element, $A1_{xx}$ and $A2_{xx}$ are the $xx^{th}$ component of tensor $\bs{A}$ in element $i=1$ and $i=2$, respectively, and $A1_{xy}$ and $A2_{xy}$ are the $xy^{th}$ component of tensor $\bs{A}$ in element $i=1$ and $i=2$, respectively. In the limit of $\Delta x \rightarrow 0$, Eq.\@ (\ref{eq:Ae5}) remains constant as
\beq
\int_{\Omega_{5}} ||\nabla\cdot\bs{A}(\bx)||_2 \ d\Omega_{5} = l \ \sqrt{\left( A2_{xx} - A1_{xx} \right)^2+\left( A2_{xy} - A1_{xy}\right)^2} \ .
\label{eq:Ae6}
\eeq
A similar calculation follows for elements $6$ through $8$ as $\Delta x \rightarrow 0$ and $\Delta y \rightarrow 0$. To simplify the integral of Eq.\@ (\ref{eq:mer}) within element $(9)$ we assume that for this regular mesh $(\nabla\cdot\bs{A}(\bx))\cdot(\nabla \cdot \bs{A}(\bx))$ can be bounded by the following expression:
\beq
||\nabla\cdot\bs{A}(\bx)||_2 \leq \sqrt{\left(\frac{\Delta_x A_{xx}}{\Delta x}+\frac{\Delta_y A_{xy}}{\Delta y}\right)^2+\left(\frac{\Delta_x A_{xy}}{\Delta x}+\frac{\Delta_y A_{yy}}{\Delta y}\right)^2}
\eeq
where
\begin{alignat}{4}
\Delta_x A_{xx} &= \text{max}(|A2_{xx}-A1_{xx}|,|A4_{xx}-A3_{xx}|),\\
\Delta_y A_{xy} &= \text{max}(|A3_{xy}-A1_{xy}|,|A4_{xy}-A2_{xy}|),\\
\Delta_x A_{xy} &= \text{max}(|A2_{xy}-A1_{xy}|,|A4_{xy}-A3_{xy}|), \text{and}\\
\Delta_y A_{yy} &= \text{max}(|A3_{yy}-A1_{yy}|,|A4_{yy}-A2_{yy}|)
\end{alignat}
are constants within the element. Thus, the integral is also bounded as follows: 
\beq
\int_{\Omega_{9}} ||\nabla\cdot\bs{A}(\bx)||_2 \ d\Omega_{9} \leq (\Delta y \ \Delta x)  \sqrt{\left(\frac{\Delta_x A_{xx}}{\Delta x}+\frac{\Delta_y A_{xy}}{\Delta y}\right)^2+\left(\frac{\Delta_x A_{xy}}{\Delta x}+\frac{\Delta_y A_{yy}}{\Delta y}\right)^2}.
\eeq
However, note that in the limit of $\Delta x \rightarrow 0$ and $\Delta y \rightarrow 0$, the right-hand side of this inequality will become increasingly small, as will the integral of element $i=9$. Thus, the center-edge elements ($i=5$ to $i=8$) are the only non-zero integrals that remain in the proposed limit. Thus, for a regular mesh with piecewise constant $\bs{A}(\bx)$, the integral in Eq.\@ (\ref{eq:mer}) is simply an evaluation of the magnitude of the jumps in $\bs{A}(\bx)$, which can be found by evaluating the strains within neighboring elements at element edges, given the nodal values and shape functions of those elements, times the length of the edge and can be written: 
\beq
    \R_{Pm}[\bu(\bx)] = \alpha_{Pm} \sum_{j=1}^{N_{e}} l^j \ ||\llbracket \bs{A}^j \rrbracket \cdot \bs{n}^j|| = \alpha_{Pm} \sum_{j=1}^{N_{e}} l^j \ \sqrt{(\llbracket \A^j[\bu] \rrbracket \cdot \mathbf{n}^j) \cdot (\llbracket \A^j[\bu] \rrbracket \cdot \mathbf{n}^j))}. \label{eq:Ae7}
\eeq
Here $j$ is an index (not a power) to all the internal edges of the domain, $N_{e}$ is the number of interior element edges, $\bs{n}^j$ is a unit vector normal to the $j^{th}$ edge, $\llbracket \bs{A}^j \rrbracket$ is the difference in $\bs{A}$ between the elements that share the $j^{th}$ edge and $l^j$ is the length of the element edge. In Eq.\@ (\ref{eq:Ae7}), the multiplication of $\llbracket \bs{A}^j \rrbracket$ with $\bs{n}^j$ acts to eliminate the partial derivatives of $\nabla\cdot\bs{A}$ along the edge of the element, which should not vary for piecewise constant $\bs{A}$ (e.g., $\partial A_{xy}/\partial y \ \text{and} \ \partial A_{yy}/\partial y$ for element edges aligned with the $y$ axis and $\partial A_{xy}/\partial x \ \text{and} \ \partial A_{xx}/\partial x$ for element edges aligned with the $x$ axis). For an irregular mesh that is not aligned with the $x$ and $y$ axis, it can be shown that Eq.\@ (\ref{eq:Ae7}) remains true provided that $\bs{A}(\bx)$ is piecewise constant. However, we do not present the proof here.

To avoid discontinuity in the derivative, we approximate Eq.\@ (\ref{eq:Ae7}) (as in Eq.\@ (\ref{eq:mreg})) as follows:
\beq
\R_{Pm}^{\delta}[\bu(\bx)] \approx \alpha_{Pm} \sum_{j=1}^{N_{e}} l^j \sqrt{(\llbracket \A^j[\bu] \rrbracket \cdot \mathbf{n}^j) \cdot (\llbracket \A^j[\bu] \rrbracket \cdot \mathbf{n}^j)) + \delta} \, , \label{eq:dmreg}
\eeq 
Although the above simplification assumes a piecewise constant $\bs{A}(\bx)$, the use of bilinear, quadrilateral shape functions in the FE discretization of $\bu(\bx)$, chosen to exploit a reduced integration scheme and avoid mesh locking, results in a function $\bs{A}(\bx)$ that is not strictly piecewise constant. Thus, the values of $\bs{A}(\bx)$ for elements adjacent to the $j^{th}$ edge were evaluated in the center of the $j^{th}$ edge and assumed to be equal to the approximate value of $\bs{A}(\bx)$ within the respective element. For the plane strain SPREME regularization, $\R_{P\epsilon}$, the value of the dilatation term of $\bs{A}(\bx)$ (i.e., $\bar{\lambda} \nabla \cdot \bu(\bx) \I$) is evaluated at the element center, rather than the edge center, to maintain a reduce integration and avoid mesh locking.

The optimization of Eq.\@ (\ref{eq:dmreg}) was performed by first creating two $N_e \times 2N_n$ matrix operators $K_{ji}^x$ and $K_{ji}^y$ which, for a given mesh with $N_n$ nodes and $N_e$ interior edges, calculate the $x$ and $y$ components of $\llbracket \A^j[\bu] \rrbracket \cdot \mathbf{n}^j$ at the center of the $j^{th}$ edge of the domain, respectively. In addition, we define an array $\overline{u}_i$ of size $2N_n$ as the nodal values of the discretized function $\bu(\bx)$ such that $\llbracket \A^j[\bu] \rrbracket \cdot \mathbf{n}^j=(K_{ji}^x \ \overline{u}_i) \bs{e}_x + (K_{ji}^y \ \overline{u}_i) \bs{e}_y$, where $\bs{e}_x$ and $\bs{e}_y$ are unit vectors in the $x$ and $y$ directions, respectively, and the summation over $i$ is implied. Then Eq.\@ (\ref{eq:dmreg}) can be written as:
\beq
\overline{\R}_{Pm}^{\delta} = \alpha_{Pm} \sum_{j=1}^{N_{e}} l^j \sqrt{(K_{ji}^x \ \overline{u}_i)^2 + (K_{ji}^y \ \overline{u}_i)^2 + \delta}, \label{eq:dmreg2}
\eeq
where $(\cdot)^2$ is performed on the each $j^{th}$ element edge separately, with no implied summation over $j$ and $\overline{\R}_{Pm}^{\delta}$ is a discrete approximation of $\R_{Pm}^{\delta}[\bu(\bx)]$. The first derivative of Eq.\@ (\ref{eq:dmreg2}) is then:
\beq
\frac{\partial \overline{\R}_{Pm}^{\delta}}{\partial \overline{u}_n} = \alpha_{Pm} \sum_{j=1}^{N_{e}} l^j \frac{(K_{ji}^x \ \overline{u}_i)K_{jn}^x+(K_{ji}^y \ \overline{u}_i)K_{jn}^y}{\sqrt{(K_{ji}^x \ \overline{u}_i)^2 + (K_{ji}^y \ \overline{u}_i)^2 + \delta}}, \label{eq:dmreg3}
\eeq
where $n$ is the index of the derivative of $\overline{\R}_{Pm}^{\delta}$ w.r.t. the $n^{th}$ value of $\overline{u}_n$. In Eq.\@ (\ref{eq:dmreg3}), there is still no implied summation of $j$ in the argument to the explicit summation of $j$, however, the implied summation over $i$ remains. The second derivative of Eq.\@ (\ref{eq:dmreg2}) is approximated as:
\beq
\frac{\partial^2 \overline{\R}_{Pm}^{\delta}}{\partial \overline{u}_n\partial \overline{u}_m} \approx \alpha_{Pm} \sum_{j=1}^{N_{e}} l^j \frac{K_{jm}^x K_{jn}^x+K_{jm}^y K_{jn}^y}{\sqrt{(K_{ji}^x \ \overline{u}_i)^2 + (K_{ji}^y \ \overline{u}_i)^2 + \delta}}, \label{eq:dmreg4}
\eeq
where $n$ and $m$ index derivatives of $\overline{\R}_{Pm}^{\delta}$ w.r.t. $\overline{u}_n$ and $\overline{u}_m$, respectively. Again, there is no implied sum of $j$ in the argument to the explicit sum of $j$, however, there is an implied summation over $i$ in the denominator. It should be noted that, similar to the continuous case in Eq.\@ (\ref{eq:gmr}), we have neglected the nonlinear term in the second derivative in Eq.\@ (\ref{eq:dmreg4}).

\subsection{Algorithmic Parameter Justification}
\label{apss:par}
This section is intended to justify our choice of regularization parameters $\alpha_{\epsilon}, \alpha_{\epsilon i}, \alpha_{P\epsilon}$, and $\alpha_{P\sigma}$ for all registrations used in this work. To do this, a single-frame registration measurement was performed between frames 1 and 7. The average axial compressive strain between these frames was approximately $0.8\%$, which corresponds to the maximum frame-to-frame strain for the simulated sequences used in this work. The resulting total strain measurement error was calculated for all types of regularization and inclusion magnitudes for a range of $\alpha$ values ($10^{\text{-}6}$ to $10^6$). All strain error values are plotted in Figure \ref{fig:alphachoice}. For each type of regularization, the value of the regularization parameter chosen for use in this paper was the average value $\alpha$ that minimized the error in all types of inclusions. The parameters were found to be $\alpha_{\epsilon}=24.8$, $\alpha_{\epsilon i}=1.89e4$, $\alpha_{P\epsilon}=1.14e\text{-}4$, and $\alpha_{P\sigma}=2.33e\text{-}4$. These regularization parameters were used for all subsequent registrations in this work with their respective regularization types.
\ifnum \Plots>0
\begin{figure}[H]
    \includegraphics[trim={0cm 0 0cm 0},clip,width = \textwidth]{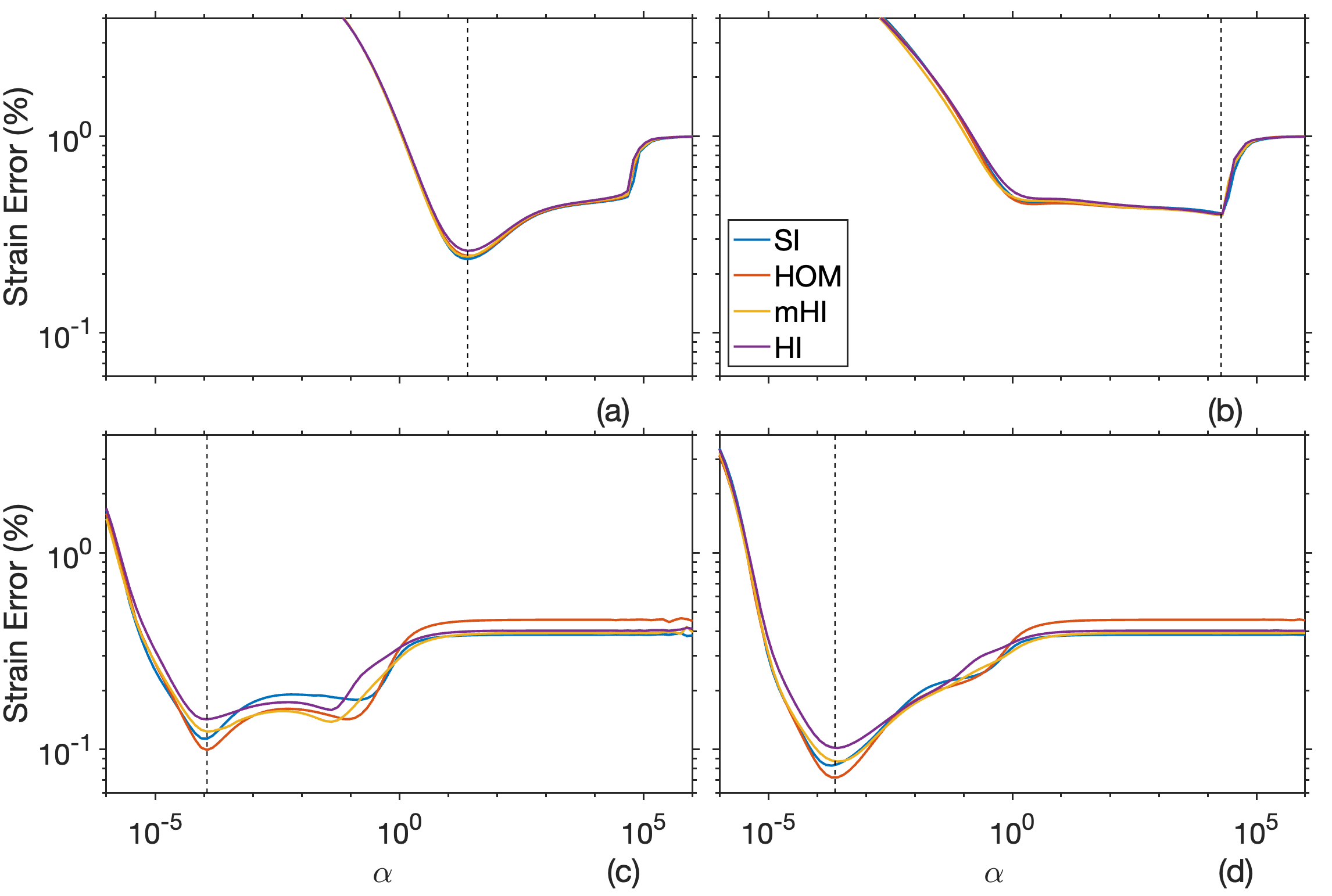}
    \caption{Plots of strain error as a function of varied regularization parameter ($\alpha$) for the 3D image set (frames 1 to 7) with all inclusion types and regularization (a) $\R_\epsilon$, (b) $\R_{\epsilon i}$, (c) $\R_{P\epsilon}$, and (d) $\R_{P\sigma}$.}
    \label{fig:alphachoice}
\end{figure}
\fi

\subsection{Supplemental Results}
\label{apss:suppres}
Additional results supporting the main findings are provided here. Table \ref{table:Astrs} shows the calculated values of the error and contrast metrics for the soft (SI), homogeneous (HOM) and medium hard inclusions (mHI) for the 3D simulated model sequences. Images of lateral and shear strains recovered from the plane strain and plane stress forward simulations (HI) are shown in Figures \ref{fig:Alat} and \ref{fig:Ashear}, respectively.

\begin{table}
\begin{center}
\begin{tabular}{|R{0.75cm}||C{\cw}|C{\cw}|C{\cw}|C{\cw}|C{\cw}|C{\cw}|C{\cw}|C{0.7cm}|C{0.95cm}|} 
 \hline
 \multicolumn{8}{|c|}{\bf{Measurement Error (Single-Frame Meas. 1 to 7, 3D Model)}}&\multicolumn{2}{|c|}{\bf{Contrast}}\\
 \hline
 \multicolumn{1}{|c||}{Model/Reg.}&\multicolumn{2}{|c|}{Disp. Comp.} &Total Disp.&\multicolumn{3}{|c|}{Strain Component} &Total Strain&\multicolumn{2}{|c|}{$\epsilon_{xx}$ Only}\\
 \hline\hline
 \multicolumn{1}{|l||}{\bf{SI}}  & $u_{x}$ & $u_{y}$  & $\mathbf{u}$ & $\epsilon_{xx}$ & $\epsilon_{yy}$ & $\epsilon_{xy}$  & $\epsilon$&$SR$&$CNR_e$\\
 \hline\hline
 $R_\epsilon$ & $0.2\%$ & $18.6\%$ & $2.7\%$ & $13.3\%$ & $41.4\%$ & $66.6\%$ & $23.8\%$ & $0.60$ & $7.15$ \\
\hline
$R_{\epsilon i}$ & $0.2\%$ & $98.0\%$ & $14.2\%$ & $4.7\%$ & $94.7\%$ & $70.5\%$ & $43.8\%$ & $0.61$ & $8.42$ \\
\hline
$R_{P\epsilon}$ & $0.1\%$ & $15.5\%$ & $2.3\%$ & $5.2\%$ & $17.4\%$ & $50.3\%$ & $11.3\%$ & $0.60$ & $9.04$ \\
\hline
$R_{P\sigma}$ & $0.1\%$ & $7.5\%$ & $1.1\%$ & $4.2\%$ & $12.9\%$ & $34.1\%$ & $8.3\%$ & $0.60$ & $9.25$ \\
 \hline\hline
 \multicolumn{1}{|l||}{\bf{HOM}} & $u_{x}$ & $u_{y}$  & $\mathbf{u}$ & $\epsilon_{xx}$ & $\epsilon_{yy}$ & $\epsilon_{xy}$  & $\epsilon$&$SR$&$CNR_e$\\
 \hline\hline
$R_\epsilon$ & $0.2\%$ & $18.3\%$ & $2.7\%$ & $13.7\%$ & $43.6\%$ & $63.0\%$ & $24.7\%$ & $0.76$ & $1.85$ \\
\hline
$R_{\epsilon i}$ & $0.2\%$ & $96.3\%$ & $13.9\%$ & $4.6\%$ & $93.6\%$ & $55.0\%$ & $43.0\%$ & $0.76$ & $2.55$ \\
\hline
$R_{P\epsilon}$ & $0.1\%$ & $13.9\%$ & $2.0\%$ & $4.3\%$ & $16.3\%$ & $40.0\%$ & $10.0\%$ & $0.76$ & $2.47$ \\
\hline
$R_{P\sigma}$ & $0.1\%$ & $5.4\%$ & $0.8\%$ & $3.4\%$ & $12.2\%$ & $25.4\%$ & $7.2\%$ & $0.77$ & $2.53$ \\
 \hline\hline
 \multicolumn{1}{|l||}{\bf{mHI}} & $u_{x}$ & $u_{y}$  & $\mathbf{u}$ & $\epsilon_{xx}$ & $\epsilon_{yy}$ & $\epsilon_{xy}$  & $\epsilon$&$SR$&$CNR_e$\\
 \hline\hline
 $R_\epsilon$ & $0.2\%$ & $19.1\%$ & $2.8\%$ & $13.9\%$ & $42.5\%$ & $59.9\%$ & $24.5\%$ & $1.08$ & $0.12$ \\
\hline
$R_{\epsilon i}$ & $0.2\%$ & $96.0\%$ & $13.9\%$ & $5.0\%$ & $92.7\%$ & $59.9\%$ & $43.0\%$ & $1.05$ & $0.05$ \\
\hline
$R_{P\epsilon}$ & $0.1\%$ & $16.3\%$ & $2.4\%$ & $5.2\%$ & $19.5\%$ & $48.8\%$ & $12.4\%$ & $1.08$ & $0.14$ \\
\hline
$R_{P\sigma}$ & $0.1\%$ & $7.3\%$ & $1.1\%$ & $4.1\%$ & $14.9\%$ & $28.5\%$ & $8.8\%$ & $1.08$ & $0.15$ \\
 \hline
\end{tabular}
\caption{Percent displacement and strain error calculated for a single-frame measurement (Frames 1 to 7) and the 3D model forward simulation with inclusion modulus values  of E = 20 MPa (mHI), E = 10 MPa (HOM), and E = 5 MPa (SI) and shown for all regularization types. Results are given for individual displacement vector and strain tensor components and total displacement and strain values.Contrast recovery is also reported as $SR$ and $CNR_e$. The strain ratio values calculated from the true strains for the forward models of the SI, HOM and mHI inclusions are $SR = 0.60$, $0.76$, and $1.07$ respectively. The corresponding contrast to noise ratios are $CNR_e = 9.64$, $2.59$, and $0.13$, respectively.}
\label{table:Astrs}
\end{center}
\end{table}

\ifnum \Plots>0
\begin{figure}[H]
    \begin{tikzpicture}
        \node (image) at (0,0) {\includegraphics[trim={1cm 0 5cm 0},clip,width = \textwidth]{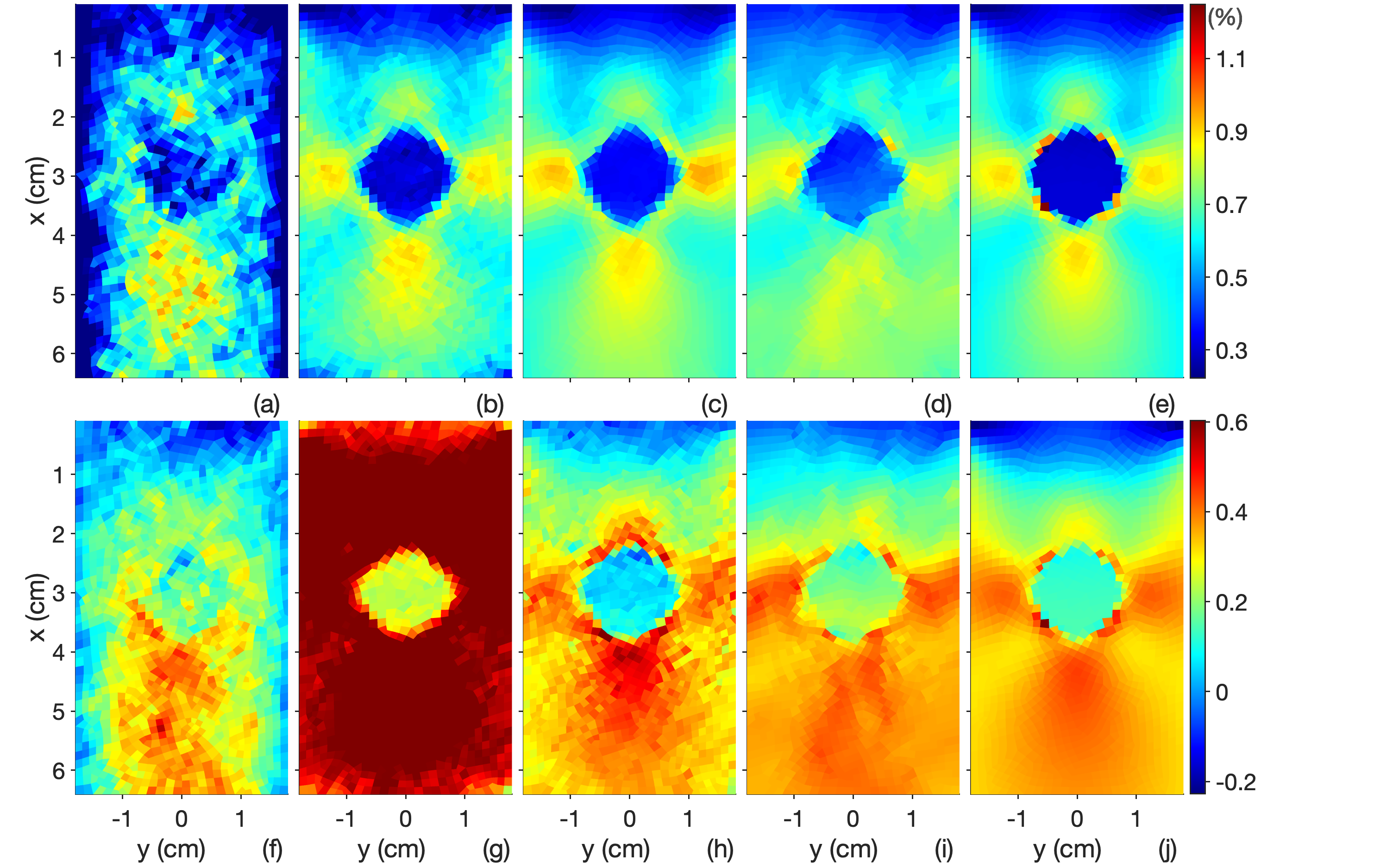}};
        \node[above, left] at (-5.1, 5.5) {$\mathcal{R}_\epsilon$};
	\node[above, left] at (-2.5, 5.5) {$\mathcal{R}_{\epsilon i}$};
	\node[above, left] at (0.3, 5.5) {$\mathcal{R}_{P\epsilon}$};
	\node[above, left] at (3, 5.5) {$\mathcal{R}_{P\sigma}$};
	\node[above, left] at (5.8, 5.5) {True};	
	\node[above, rotate=90] at (-7.5, 2.75) {$\epsilon_{yy}$ (Fwd. Plane Strain)};
	\node[above, rotate=90] at (-7.5, -2.5) {$\epsilon_{yy}$ (Fwd. Plane Stress)};
    \end{tikzpicture}        
    \caption{(a-d) Lateral strain ($\epsilon_{yy}$) measurements for 2D plane strain image set (frames 1 to 7, HI) and regularization (a) $\R_\epsilon$, (b) $\R_{\epsilon i}$, (c) $\R_{P\epsilon}$, and (d) $\R_{P\sigma}$. (e) simulated lateral strain for 2D plane strain. (f-i) Lateral strain ($\epsilon_{yy}$) measurements for 2D plane stress image set (frames 1 to 7, HI) and regularization (f) $\R_\epsilon$, (g) $\R_{\epsilon i}$, (h) $\R_{P\epsilon}$, (i) $\R_{P\sigma}$, and (j) simulated lateral strain for 2D plane stress.}
    \label{fig:Alat}
\end{figure}
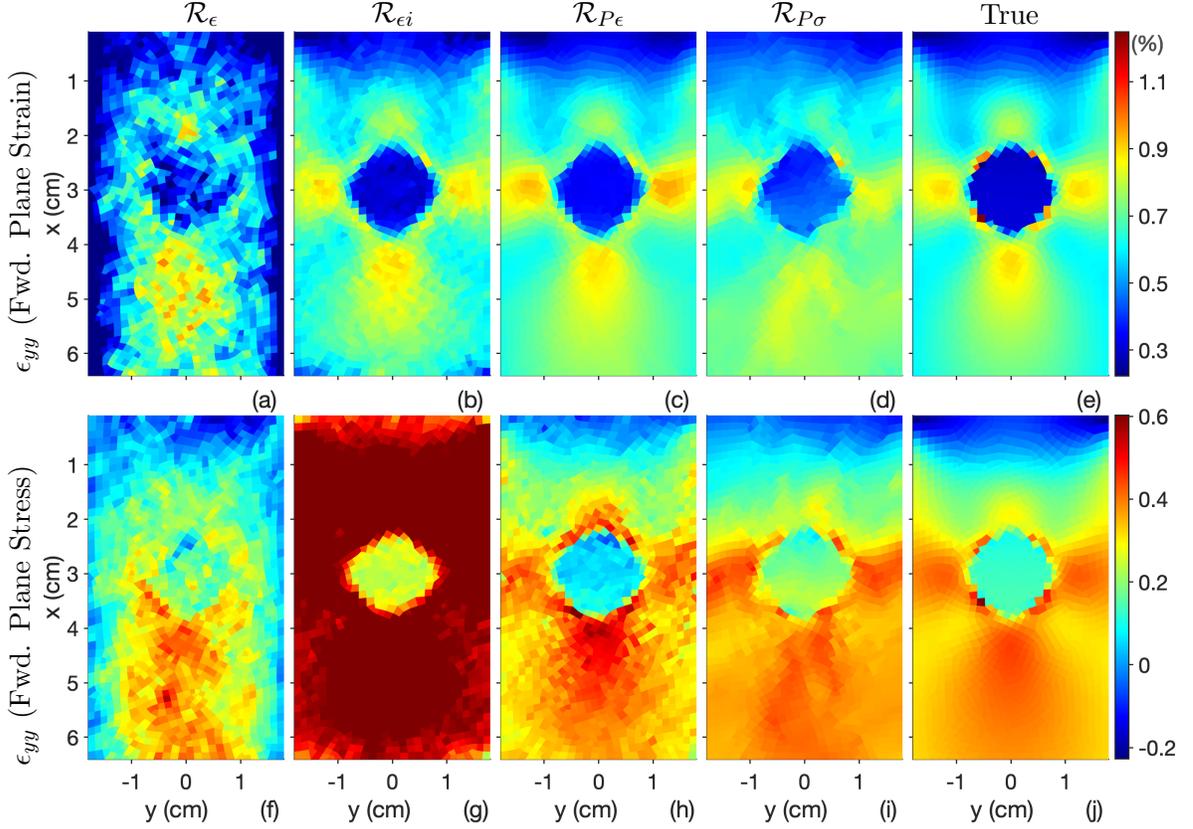
\fi

\ifnum \Plots>0
\begin{figure}[H]
    \begin{tikzpicture}
        \node (image) at (0,0) {\includegraphics[trim={1cm 0 5cm 0},clip,width = \textwidth]{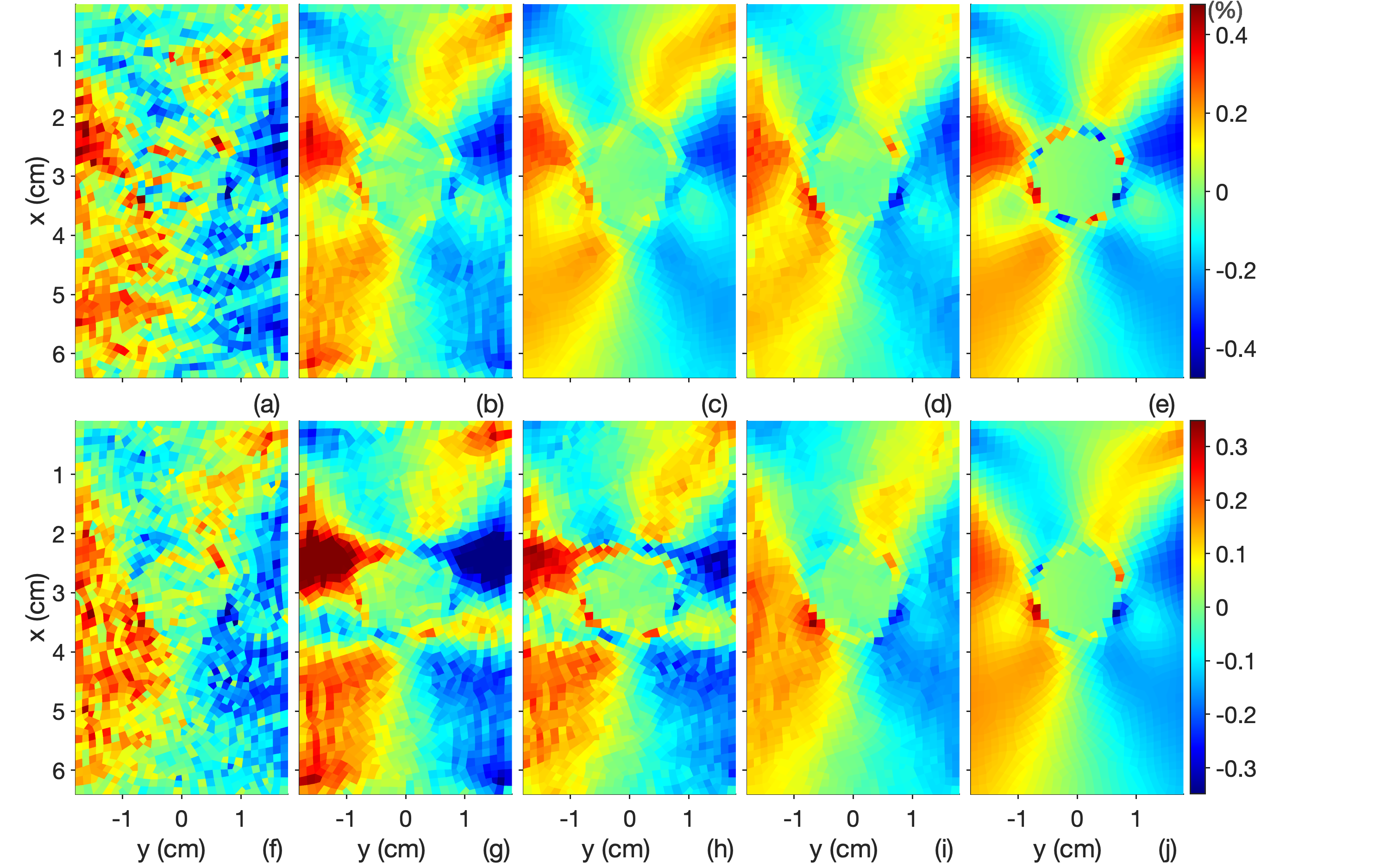}};
        \node[above, left] at (-5.1, 5.5) {$\mathcal{R}_\epsilon$};
	\node[above, left] at (-2.5, 5.5) {$\mathcal{R}_{\epsilon i}$};
	\node[above, left] at (0.3, 5.5) {$\mathcal{R}_{P\epsilon}$};
	\node[above, left] at (3, 5.5) {$\mathcal{R}_{P\sigma}$};
	\node[above, left] at (5.8, 5.5) {True};	
	\node[above, rotate=90] at (-7.5, 2.75) {$\epsilon_{xy}$ (Fwd. Plane Strain)};
	\node[above, rotate=90] at (-7.5, -2.5) {$\epsilon_{xy}$ (Fwd. Plane Stress)};
    \end{tikzpicture}
    \caption{(a-d) Shear strain ($\epsilon_{yy}$) measurements for 2D plane strain image set (frames 1 to 7, HI) and regularization (a) $\R_\epsilon$, (b) $\R_{\epsilon i}$, (c) $\R_{P\epsilon}$, (d) $\R_{P\sigma}$, (e) simulated shear strain for 2D plane strain. (f-i) Shear Strain plots for 2D plane stress image set (frames 1 to 7, HI) and regularization (f) $\R_\epsilon$, (g) $\R_{\epsilon i}$, (h) $\R_{P\epsilon}$, (i) $\R_{P\sigma}$, (j) simulated shear strain for 2D plane stress.}
    \label{fig:Ashear}
\end{figure}
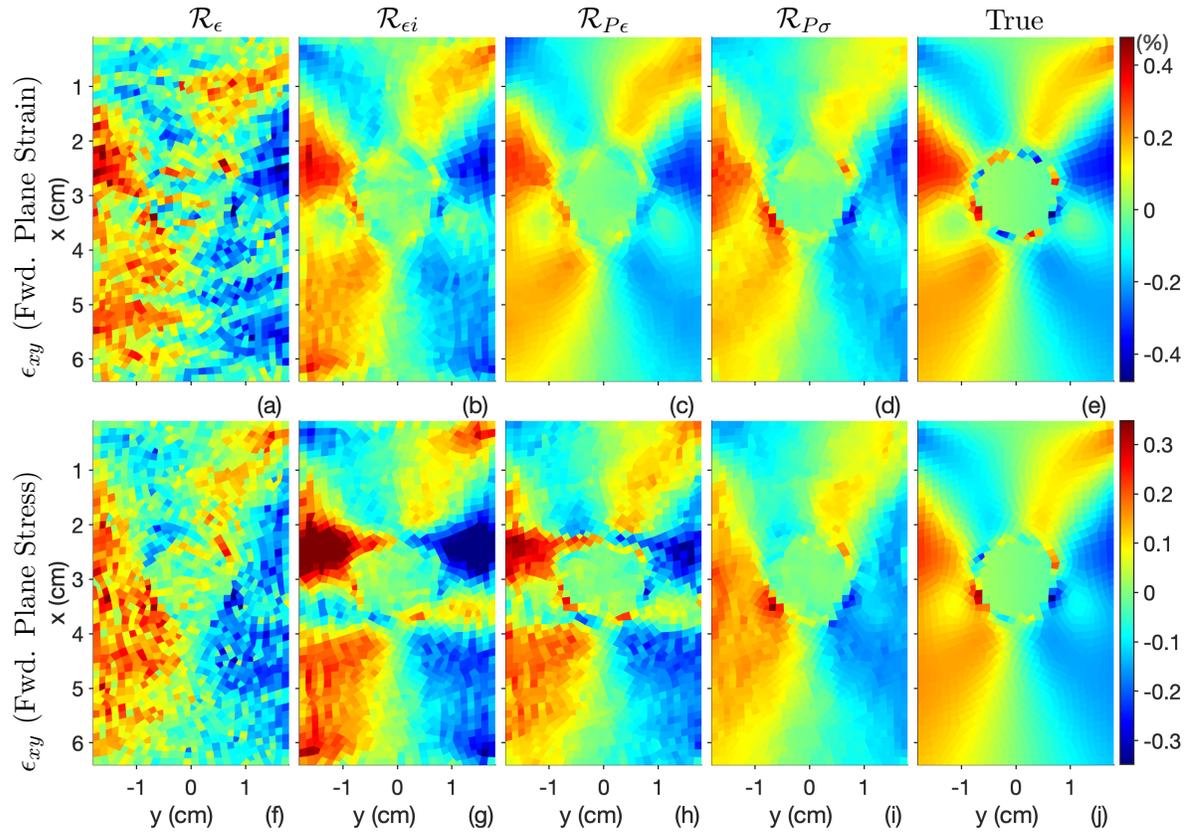
\fi

\clearpage
\bibliographystyle{apalike}
\bibliography{spremeIR_paper}
\end{document}